\documentclass[twoside,10pt]{article}
\usepackage{amssymb,amsmath,euscript}
\usepackage{graphicx,color,caption}

\usepackage{algorithm,algpseudocode}

\newtheorem{proposition}{Proposition}[section]
\newtheorem{theorem}[proposition]{Theorem}

\newtheorem{corollary}[proposition]{Corollary}
\newtheorem{definition}[proposition]{Definition}
\newtheorem{remark}[proposition]{Remark}
\newtheorem{example}[proposition]{Example}

\newcommand{\qed}{\hphantom{.}\hfill $\Box$\medbreak}

\def\S{\mathbb{S}}

\def\R{\mathbb{R}}

\def\T{\mathbb{T}}
\def\Q{\mathcal{Q}}
\def\A{\mathcal{A}}

\def\B{\mathcal{B}}
\def\C{\mathcal{C}}

\def\I{\mathcal{I}}
\title{\bf{Doubly Nonnegative Tensors, Completely Positive Tensors and Applications}
\thanks{This research was supported by the National Natural Science Foundation of China (11301022,11431002), the State Key Laboratory of Rail Traffic Control and Safety, Beijing Jiaotong University (RCS2014ZT20, RCS2014ZZ01), and the Hong Kong Research Grant Council (Grant No.
PolyU 502111, 501212, 501913 and 15302114).}}
\author{Ziyan
Luo,  \thanks{ State Key Laboratory of Rail Traffic Control and Safety, Beijing
Jiaotong University, Beijing 100044, P.R. China; ({\tt
starkeynature@hotmail.com}).} \hspace{1mm}  Liqun Qi \thanks{Department of Applied Mathematics, The Hong Kong Polytechnic University, Hung Hom, Kowloon, Hong Kong. ({\tt
liqun.qi@polyu.edu.hk}).}
}

\begin{document}
\maketitle

\begin{abstract}
The concept of double nonnegativity of matrices is generalized to doubly nonnegative tensors by means of the nonnegativity of all entries and $H$-eigenvalues.
This generalization is defined for tensors of any order (even or odd), while it reduces to the class of nonnegative positive semidefinite tensors in the even order case.
We show that many nonnegative structured tensors, which are positive semidefinite in the even order case, are indeed doubly nonnegative as well in the odd order case.
As an important subclass of doubly nonnegative tensors, the completely positive tensors are further studied. By using dominance properties for completely positive tensors,
we can easily exclude some doubly nonnegative tensors, such as the signless Laplacian tensor of a nonempty $m$-uniform hypergraph with $m\geq 3$, from the class of completely
positive tensors. Properties of the doubly nonnegative tensor cone and the completely positive tensor cone are established. Their relation and difference are discussed.
These show us a different phenomenon comparing to the matrix case. By employing the proposed properties, more subclasses of these two types of tensors are identified.
Particularly, all positive Cauchy tensors with any order are shown to be completely positive. This gives an easily constructible subclass of completely positive tensors,
which is significant for the study of completely positive tensor decomposition. A preprocessed Fan-Zhou algorithm is proposed which can efficiently verify the complete positivity
of nonnegative symmetric tensors.  We also give the solution analysis of tensor complementarity problems with the strongly doubly nonnegative tensor structure.


\end{abstract}


\section{Introduction}
Doubly nonnegative matrices and completely positive matrices have attracted considerable attention due to their applications in optimization, especially in creating convex formulations of NP-hard problems, such as the quadratic assignment problem in combinatorial optimization and the polynomial optimization problems (see \cite{AKK2013,AKKT2014-1,AKKT2014-2,KKT2013,YM2010,ZF2014} and references therein). In recent years, an emerging interest in the assets of multilinear algebra has been concentrated on the higher-order tensors, which serve as a numerical tool, complementary to the arsenal of existing matrix techniques. In this vein, the concept of completely positive matrices has been extended for higher-order tensors and properties of this special nonnegative tensors have been exploited in \cite{QXX2014}. In this paper, the double nonnegativity will be extended from matrices to higher-order tensors, which contains the completely positive tensors as an important subclass. Properties on doubly nonnegative tensors and completely positive tensors will be exploited and further developed.

It is well-known that the double nonnegativity of a symmetric matrix refers to the nonnegativity of all entries and eigenvalues \cite{BP1994,BM2003}. Thus, it is reasonable to define the doubly nonnegative tensor by the nonnegativity of all entries and $H$-eigenvalues of a symmetric tensor. Formally, a doubly nonnegative tensor is defined as follows.
\begin{definition}\label{def1} A symmetric tensor is said to be a \emph{doubly nonnegative} tensor if all of its entries and $H$-eigenvalues are nonnegative; A symmetric tensor is said to be a \emph{strongly doubly nonnegative} tensor if all of its entries are nonnegative and all of its $H$-eigenvalues are positive.
\end{definition}

Noting that nonnegative tensors always have $H$-eigenvalues \cite{YY2010}, the above concepts are well-defined for tensors of any order (even or odd). Particularly, inspired by the criterion for positive semidefinite (definite) tensors as shown in \cite{Q2005}, a (strongly) doubly nonnegative tensor is exactly a positive semidefinite (definite) tensor with nonnegative entries when $m$ is even. Thus, many even-order symmetric tensors, such as the diagonally dominant tensors, the generalized diagonally dominant tensors, the $H$-tensors, the complete Hankel tensors, the strong Hankel tensors, the $MB_0$-tensors, the quasi-double $B_0$-tensors, the double $B_0$-tensors, and the $B_0$-tensors, which are positive semidefinite, are therefore doubly nonnegative if they are symmetric and nonnegative.  A natural question is: Are the aforementioned structured positive semidefinite tensors with nonnegative entries still doubly nonnegative in the odd order case? An affirmative answer will be given in this paper.

Odd order tensors are extensively encountered in many fields, such as the odd-uniform hypergraphs, the diffusion tensor imaging, and signal processing. As an important concept in spectral hypergraph theory, the largest Laplacian $H$-eigenvalue and its relation to the largest signless Laplacian $H$-eigenvalue for odd-uniform hypergraphs have been studied and characterized \cite{HQS2013,HQX2015,Qi2014}, which shows a different phenomenon from that in the even case. In non-Gaussian diffusion tensor imaging, the even-order tensors only affect the magnitude of the signal, whereas odd-order tensors affect the phase of the signal \cite{LBAM2004}. And the independent component analysis, a computational method for separating a multivariate signal into additive subcomponents in signal processing, adopts a third-order tensor diagonalization which is shown to be very efficient since for third-order tensors, the computation of an elementary Jacobi-rotation is again equivalent to the best rank-one approximation just like  the case of fourth order tensors does \cite{LDV2001}. All these wide applications indicate the significance and merits of odd order tensors in tensor theory and analysis. As a fundamental and essential ingredient of tensor analysis, positive semidefiniteness has been extensively studied for even order tensors \cite{CQ2015,DQW2015,KMB2015,LL2015,LQL2015,LWZZL2014,LQY2015,Q2005,Qi2014,Q2015,QS2014,QY2014,ZQZ2014}, but it vanishes for odd order tensors. The concept of doubly nonnegative tensors as proposed in Definition \ref{def1}, to some extent, has partially made up this deficiency in theory.

Analogous to the matrix case, completely positive tensors form an extremely important part of doubly nonnegative tensors. They are connected with nonnegative tensor factorization and have wide applications \cite{CZPA, FZ2014, QXX2014, SH}.    As an extension of the completely positive matrix, a completely positive tensor admits its definition in a pretty natural way as initiated by Qi et al. in \cite{QXX2014} and recalled below.
\begin{definition} \label{def2} A tensor $\A\in\S_{m,n}$ is said to be a \emph{completely positive} tensor if there exist an integer $r$ and some $u^{(k)}\in\R^n_+$, $k\in[r]$ such that $\A=\sum\limits_{k=1}^r \left(u^{(k)}\right)^m.$ A tensor $\A\in\S_{m,n}$ is said to be a \emph{strongly completely positive} tensor if there exist an integer $r$ and some $u^{(k)}\in\R^n_+$, $k\in[r]$ such that $\A=\sum\limits_{k=1}^r \left(u^{(k)}\right)^m$ and $span\{u^{(1)},\ldots, u^{(r)}\}=\R^n$.  \end{definition}
It is worth pointing out that both in Definitions \ref{def1} and \ref{def2}, the second parts, which are new concepts developed in this paper, can be regarded as counterparts tailored for all order tensors to the positive definiteness solely customized for even order tensors.

All completely positive tensors contribute a closed convex cone in symmetric tensor space, associated with the copositive tensor cone as its dual \cite{PVZ2014,QXX2014}. Optimization programming over these two closed convex cones in the matrix case has been widely studied in the community of combinatorial optimization and quadratic programming \cite{AKKT2014-1,B2012,D2010,HS2010}. In spite of the better tightness of completely positive cone relaxation comparing to the well-known positive semidefinite relaxation, the former one is not computationally tractable. As a popular relaxation strategy, the doubly nonnegative matrix cone is always treated as a surrogate to the completely positive matrix cone due to its tractability of the involved double nonnegativity \cite{YM2010}. Many research works on these two cones emerged, not only on the algebraic and geometric properties of the cones \cite{Dickinson2010,Dickinson2011,DS2008,GST2013}, but also the relation and the difference between them \cite{BAD2009,DA2013,SZ1993}. Along this line of research, the completely positive tensor cone was employed to reformulate polynomial optimization problems which are not necessarily quadratic \cite{PVZ2014}. Numerical optimization for the best fit of completely positive tensors with given length of decomposition was formulated as a nonnegative constrained least-squares problem in Kolda's paper \cite{K2015}. A verification approach in terms of truncated moment sequences for checking completely positive tensors was proposed and an optimization algorithm based on semidefinite relaxation for completely positive tensor decomposition was established by Fan and Zhou in their recent work \cite{FZ2014}. In this paper, the properties, together with their relation and difference will also be investigated for the doubly nonnegative tensor cone and the completely positive tensor cone, which somehow shows a different phenomenon from the matrix case. A noteworthy observation is the dominance properties of completely positive tensors. These properties turn out to be a very powerful tool to exclude some higher-order tensors, such as the well-known signless Laplacian tensors of nonempty $m$-uniform hypergraphs with $m\geq 3$, from the class of completely positive tensors. More subclasses of doubly nonnegative tensors and completely positive tensors will be discussed. In particular, we show that positive Cauchy tensors of any order are completely positive. This provides another sufficient condition for completely positive tensors. As an application, a preprocessing scheme for checking complete positive tensors based on the zero-entry dominance property and a simplified strong dominance property is proposed. Another application is the solution analysis of tensor complementarity problems with strongly doubly nonnegative tensors.

The rest of the paper is organized as follows. In Section 2, we briefly review some basic concepts and properties on symmetric tensors, and definitions of several structured tensors that will be discussed in the subsequent analysis.
The doubly nonnegative tensors and their properties will be discussed in Section 3, where we show that several structured tensors, which are positive semidefinite in the even order case, are all doubly nonnegative in the odd order
case as well, provided that they are symmetric and entry-wise nonnegative. Section 4 is devoted to the completely positive tensors and their properties, including the dominance properties which are very useful to exclude some tensors
from the complete positivity. Particularly, an interesting finding is all signless Laplician tensors of nonempty $m$-uniform hypergraphs with $m\geq 3$ are not completely positive. In order to get a further understanding on doubly
nonnegative tensors and completely positive tensors, the doubly nonnegative tensor cone and the completely positive tensor cone, together with their geometric and algebraic properties, are discussed in Section 5. The relation and
difference of these two cones are also analyzed which bring out a more complicated phenomenon compared to the matrix case. More subclasses of nonnegative tensors and completely positive tensors are investigated in Section 6.
A preprocessing scheme for completely positive tensor verification and decomposition, along with some numerical tests, is stated in Section 7. The solution existence of tensor complementarity problems with the strongly doubly nonnegative tensor structure is characterized as well. Concluding remarks are drawn in Section 8.

Some notations that will be used throughout the paper are listed here. The $n$-dimensional real Euclidean space is denoted by $\R^n$, where $n$ is a given natural number. The nonnegative orthant in $\R^n$ is denoted by $\R^n_+$, with the interior $\R^n_{++}$ consisting of all positive vectors. The $n$-by-$l$ real matrix space is denoted by $\R^{n\times l}$. Denote $[n]:=\{1,2,\ldots, n\}$. Vectors are denoted by lowercase letters such as $x$, $u$, matrices are denoted by capital letters such as $A$, $P$, and tensors are written as calligraphic capital letters such as $\A$, $\B$. The space of all real $m$th order $n$-dimensional tensors is denoted by $\T_{m,n}$, and the space of all symmetric tensors in $\T_{m,n}$ is denoted by $\S_{m,n}$. For any closed convex set $\mathbb{M}$, a closed convex subset $\mathbb{F}$ is said to be a face of $\mathbb{M}$ if any $x$, $y\in {\mathbb{M}}$ satisfying $x+y\in {\mathbb{F}}$ implies that $x$, $y\in {\mathbb{F}}$. We use ${\mathbb{F}}\unlhd {\mathbb{M}}$ if ${\mathbb{F}}$ is a face of $\mathbb{M}$. For a subset $\Gamma\subseteq [n]$, $|\Gamma|$ stands for its cardinality.


\section{Preliminaries}

\subsection{Symmetric Tensors and Eigenvalues}
Let $\A=\left(a_{i_1\ldots i_m}\right)$ be an $m$th order $n$-dimensional real tensor. $\A$ is called a \emph{symmetric} tensor if the entries $a_{i_1\ldots i_m}$ are invariant under any permutation of their indices for all $i_j\in [n]$ and $j\in [m]$, denoted as $\A\in \S_{m,n}$. A symmetric tensor $\A$ is said to be positive semidefinite (definite) if $\A x^m:=\sum\limits_{i_1,\ldots,i_m\in [n]}a_{i_1\ldots i_m}x_{i_1}\cdots x_{i_m}\geq 0 (>0)$ for any $x\in \R^n\setminus \{0\}$. Here, $x^m$ is a rank-one tensor in $\S_{m,n}$ defined as $\left(x^m\right)_{i_1\ldots i_m}:=x_{i_1}\cdots x_{i_m}$ for all $i_1$, $\ldots$, $i_m\in [n]$. Evidently, when $m$ is odd, $\A$ could not be positive definite and $\A$ is positive semidefinite if and only if $\A=O$, where $O$ stands for the zero tensor. A tensor $\A\in \T_{m,n}$ is said to be (strictly) copositive if $\A x^m\geq 0$ ($>0$) for all $x\in\R^n_{+}\setminus \{0\}$. The definitions on eigenvalues of symmetric tensors are recalled as follows.

\begin{definition}[\cite{Q2005}] Let $\A\in \S_{m,n}$ and $\mathbb{C}$ be the complex field. We say that $(\lambda,x)\in {\mathbb{C}}\times \left({\mathbb{C}}^n\setminus \{0\}\right)$ is an eigenvalue-eigenvector pair of $\A$ if
$\A x^{m-1}=\lambda x^{[m-1]},$ where $\A x^{m-1}$ and $x^{[m-1]}$ are all $n$-dimensional column vectors given by
$$\left(\A x^{m-1}\right)_i:=\sum\limits_{i_2,\ldots,i_m\in[n]} a_{ii_2\ldots i_m} x_{i_2}\cdots x_{i_m},~\left(x^{[m-1]}\right)_i=x_i^{m-1},~~\forall i\in [n].\eqno(2.1)$$  If the eigenvalue $\lambda$ and the eigenvector $x$ are real, then $\lambda$ is called an $H$-eigenvalue of $\A$ and $x$ an $H$-eigenvector of $\A$ associated with $\lambda$. If $x\in \R^n_+ (\R^n_{++})$, then $\lambda$ is called an $H^+ (H^{++})$-eigenvalue of $\A$.
\end{definition}

\begin{definition}[\cite{Q2005}] Let $\A\in \S_{m,n}$ and $\mathbb{C}$ be the complex field. We say that $(\lambda,x)\in {\mathbb{C}}\times \left({\mathbb{C}}^n\setminus \{0\}\right)$ is an $E$ eigenvalue-eigenvector pair of $\A$ if
$\A x^{m-1}=\lambda x$ and $x^{T}x=1$, where $\A x^{m-1}$ is defined as in (2.1). If the $E$-eigenvalue $\lambda$ and the eigenvector $x$ are real, then $\lambda$ is called a $Z$-eigenvalue of $\A$ and $x$ a $Z$-eigenvector of $\A$ associated with $\lambda$.
\end{definition}

\subsection{Structured Tensors}
Several structured tensors are recalled which will be studied in the sequel of the paper.

\begin{definition}[Definition 3.14, \cite{ZQZ2014}]\label{def-diagonal} We call $\A\in \T_{m,n}$ a \emph{diagonally dominant} tensor if
$$|a_{ii\ldots i}|\geq \sum\limits_{(i_2,\ldots,i_m)\neq (i,\ldots, i)}|a_{ii_2\ldots i_m}|, ~\forall i\in [n].\eqno(2.2)$$
$\A$ is said to be \emph{strictly diagonally dominant} if the strict inequality holds in (2.2) for all $i\in [n]$.
\end{definition}

For any $\A \in \S_{m,n}$, and any invertible diagonal matrix $D=diag(d_1,\ldots, d_n)$, define $\A D^{1-m} \underbrace{D\cdots D}\limits_{m-1}$ as
$$\left(\A D^{1-m} \underbrace{D\cdots D}\limits_{m-1}\right)_{i_1\ldots i_m}=a_{i_1\ldots i_m}d_{i_1}^{1-m}d_{i_2}\cdots d_{i_m}, ~\forall i_j\in [n],~j\in [m].\eqno(2.3)$$
\begin{definition}[\cite{DQW2013,KMB2015}] A tensor $\A$ is called \emph{(strictly) generalized diagonally dominant} if there exists some positive diagonal matrix $D$ such that $\A D^{1-m} \underbrace{D\cdots D}\limits_{m-1}$, as defined in (2.3), is (strictly) diagonally dominant.
\end{definition}

\begin{definition}[\cite{Q2015}]\label{def-Hankel} Let $\A=\left(a_{i_1\ldots i_m}\right) \in \T_{m,n}$. If there is a vector $v=(v_0,\ldots, v_{(n-1)m})^{T}\in\R^{(n-1)m+1}$ such that
 $$a_{i_1\ldots i_m}=v_{i_1+\cdots+i_m-m},~\forall i_j\in [n],~j\in [m],\eqno(2.4)$$
then we say that $\A$ is an $m$th order $n$-dimensional \emph{Hankel tensor}. Let $A=(a_{ij})$ be a $\lceil\frac{(n-1)m+2}{2}\rceil \times \lceil\frac{(n-1)m+2}{2}\rceil$ Hankel matrix with $a_{ij}:=v_{i+j-2}$, where $v_{2\lceil\frac{(n-1)m+2}{2}\rceil}$ is an additional number when
$(n-1)m$ is odd. If $A$ is positive semidefinite, then $\A$ is called a \emph{strong Hankel tensor}. Suppose $\A$ is a Hankel tensor with its vandermonde decomposition $\A=\sum\limits_{k=1}^r \alpha_k \left(u^{(k)}\right)^m$, where $u^{(k)}:=(1,\xi_k,\ldots, \xi_k^{n-1})^T$, $\xi_k\in\R$, for all $k\in [r]$. If $\alpha_k>0$ for all $k\in r$, then $\A$ is called a \emph{complete Hankel tensor}.
\end{definition}

\begin{definition}[\cite{DQW2013,KMB2015}] \label{def-M-H} A tensor $\A$ is called a $Z$-tensor if there exists a nonnegative tensor $\B$ and a real number $s$ such that $\A =s\I-\B$. A $Z$-tensor $\A =s\I-\B$ is said to be an \emph{$M$-tensor} if $s\geq \rho(\B)$, where $\rho(\B)$ is the spectral radius of $\B$. If $s>\rho(\B)$, then $\A$ is called a \emph{strong $M$-tensor}. The \emph{comparison tensor} of a tensor $\A=(a_{i_1\ldots i_m})\in \T_{m,n}$, denoted by $M(\A)$, is defined as
$$\left(M(\A)\right)_{i_1\ldots i_m}:=\left\{
                                        \begin{array}{ll}
                                          |a_{i_1\ldots i_m}|, & \hbox{if ~$i_1=\ldots=i_m$;} \\
                                          -|a_{i_1\ldots i_m}|, & \hbox{otherwise.}
                                        \end{array}
                                      \right.$$
$\A$ is called an $H$-tensor (strong $H$-tensor) if its comparison tensor $M(\A)$ is an $M$-tensor (strong $M$-tensor).
\end{definition}

\begin{definition}[\cite{QS2014}] \label{def-B0} A tensor $\A=\left(a_{i_1\ldots i_m}\right) \in \T_{m,n}$ is called a $B_0$-tensor if
$$ \sum\limits_{i_2,\ldots, i_m\in [n]}a_{ii_2\ldots i_m} \geq 0$$ and $$\frac{1}{n^{m-1}}\sum\limits_{i_2,\ldots,i_m\in [n]}a_{ii_2\ldots i_m}\geq a_{ij_2\ldots j_m},~\forall (j_2,\ldots,j_m)\neq (i,\ldots,i).$$
\end{definition}
\vskip 2mm

In \cite{LL2015,LQL2015}, the $B_0 (B)$-tensor has been generalized and further studied. For any tensor $\A=\left(a_{i_1\ldots i_m}\right) \in \T_{m,n}$, denote $$\beta_i(\A):=\max\limits_{
\begin{subarray}{c}
j_2,\ldots,j_m\in [n]  \\
(i,j_2,\ldots,j_m)\neq (i,i,\ldots,i)
\end{subarray}}\{0,a_{ij_2\ldots j_m}\},~\Delta_i(\A):=\sum\limits_{
\begin{subarray}{c}
j_2,\ldots,j_m\in [n]  \\
(i,j_2,\ldots,j_m)\neq (i,i,\ldots,i)
\end{subarray}}(\beta_i(\A)-a_{ij_2\ldots j_m}),$$ and $\Delta_j^i:=\Delta_j(\A)-(\beta_j(\A)-a_{jii\ldots i}),~i\neq j.$ As defined in \cite{LL2015}, $\A$ is called a \emph{double $B$-tensor} if $a_{i\ldots i}>\beta_i(\A)$ for all $i\in [n]$ and for all $i$, $j\in [n]$, $i\neq j$, $a_{i\ldots i}\geq \Delta_i(\A)$ and $(a_{i\ldots i}-\beta_i(\A))(a_{j\ldots j}-\beta_j(\A))>\Delta_i(\A)\Delta_j(\A)$. If $a_{i\ldots i}>\beta_i(\A)$ for all $i\in [n]$, and
$$(a_{i\ldots i}-\beta_i(\A))(a_{j\ldots j}-\beta_j(\A)-\Delta_j^i(\A))\geq (\beta_j(\A)-a_{ji\ldots i})\Delta_i(\A),$$ then $\A$ is called a \emph{quasi-double $B_0$ tensor}. Let $\A=\left(a_{i_1\ldots i_m}\right) \in \T_{m,n}$ and set $b_{i_1\ldots i_m}=a_{i_1\ldots i_m}-\beta_{i_1}(\A)$ for any $i_j\in [n]$ and $j\in [m]$. If $\B:=\left(b_{i_1\ldots i_m}\right)$ is an $M$-tensor, then $\A$ is called an \emph{$MB_0$-tensor}. If $\B$ is a strong $M$-tensor, then $\A$ is called an \emph{$MB$-tensor}.

\begin{definition}[\cite{CQ2015}] \label{def-cauchy} Let $c=(c_1,\ldots, c_n)^T\in\R^n$ with $c_i\neq 0$ for all $i\in [n]$. Suppose that $\C=(c_{i_1\ldots i_m})\in \T_{m,n}$ is defined as
$$c_{i_1\ldots i_m}=\frac{1}{c_{i_1}+\cdots +c_{i_m}}, ~\forall i_j\in [n], ~j\in [m].$$
Then, we say that $\C$ is an $m$th order $n$-dimensional \emph{symmetric Cauchy tensor} and the vector $c=(c_1,\ldots, c_n)^T\in\R^n$ is called the \emph{generating vector} of $\C$.
\end{definition}

\begin{definition}[\cite{CQi2015}] \label{def-circulant} Let $\A=\left(a_{i_1\ldots i_m}\right) \in \T_{m,n}$. If for all $i_j\in [n-1]$ and all $j\in [m]$,
we have $a_{i_1\ldots i_m}=a_{i_1+1\ldots i_m+1}$, then $\A$ is called an $m$th order $n$-dimensional \emph{Toeplitz tensor}. If for $i_j$, $k_j\in [n]$, $k_j=i_j+1 ~mod~(n)$, $j\in [m]$, we have  $a_{i_1\ldots i_m}=a_{k_1\ldots k_m}$, then $\A$ is called an $m$th order $n$-dimensional \emph{circulant tensor}. For a circulant tensor $\A$, let $\A_k:=\left(a^{(k)}_{j_1\ldots j_{m-1}}\right)\in \T_{m-1,n}$ be defined as $a^{(k)}_{j_1\ldots j_{m-1}}:=a_{kj_1\ldots j_{m-1}}$. $\A_k$ is said to be the \emph{$k$th row tensor} of $\A$ for $k\in [n]$. Particularly, $\A_1$ is called the root tensor of $\A$ and $c_0:=a_{1\ldots 1}$ is called the \emph{diagonal entry} of $\A$.
\end{definition}

\begin{definition}[\cite{Qi2014}] \label{def-signless} Let $G =(V, E)$ be an $m$-uniform hypergraph. The \emph{adjacency tensor} of $G$ is defined as the $m$th order $n$-dimensional tensor $\A$ whose $(i_1, \ldots, i_m)$th entry is
$$a_{i_1\ldots i_m}=\left\{
  \begin{array}{ll}
    \frac{1}{(m-1)!}, & \hbox{if~$\{i_1,\ldots, i_m\}\in E$;} \\
    0, & \hbox{otherwise.}
  \end{array}
\right.$$
Let $\mathcal{D}$ be an $m$th order $n$-dimensional diagonal tensor with its diagonal element $d_{i\ldots i}$ being $d_i$, the degree of vertex $i$, for all $i\in [n]$. Then ${\mathcal{Q}}:={\mathcal{D}}+\A$ is called the \emph{signless Laplacian tensor} of the hypergraph $G$.
\end{definition}

\section{Doubly Nonnegative Tensors}
Properties on doubly nonnegative tensors are discussed in this section with mainly two parts. The first part is devoted to some necessary and/or sufficient conditions for doubly nonnegative tensors, and the second part is dedicated to establishing the double nonnegativity of many structured tensors.

\subsection{Necessary and/or Sufficient Conditions}
We start with some necessary and sufficient conditions for doubly nonnegative tensors by means of some special types of linear operators. Besides the one defined as in (2.3), we also invoke the following linear operator as introduced in \cite{Q2005}, which possesses the decomposition invariance property as discussed in \cite{LQY2015}:
Let $P=(p_{ij})\in\R^{n\times n}$. Define a linear operator $P^m$ as
$$\left(P^m \A\right)_{i_1\ldots i_m}:=\sum\limits_{j_1,\ldots, j_m\in [n]}p_{i_1j_1}\cdots p_{i_mj_m}a_{j_1\ldots j_m},~\forall \A=\left(a_{j_1\ldots j_m}\right)\in\S_{m,n}.$$

\begin{proposition}\label{invariance} Suppose $\A\in \S_{m,n}$ is a nonnegative symmetric tensor. We have

\begin{itemize}
\item[(i)] for any given positive diagonal matrix $D\in \R^{n\times n}$, $\A$ is (strongly) doubly nonnegative if and only if $\A D^{1-m} \underbrace{D\cdots D}\limits_{m-1}$ is nonnegative and all of its $H$-eigenvalues are nonnegative (positive);

\item[(ii)] for any permutation matrix $P\in \R^{n\times n}$, $\A$ is (strongly) doubly nonnegative  if and only if $P^{m}\A$ is (strongly) doubly nonnegative.
\end{itemize}
\end{proposition}
\noindent{\bf Proof.} (i) follows directly from the fact that the eigenvalues of $\A$ coincide with those of $\A D^{1-m} \underbrace{D\cdots D}\limits_{m-1}$ as you can see in Proposition 2.1 of \cite{KMB2015}. (ii) The nonnegativity is trivial. For any permutation matrix $P$ and any $x\in \R^n$, since $\A x^{m-1}=\lambda x^{m-1}$ if and only if $P^m \A (Px)^{m-1}=\lambda (Px)^{m-1}$, one has that $x$ is an $H$-eigenvector of tensor A associated with $\lambda$ if and only if $Px$ is an $H$-eigenvector of tensor $P^m \A$ associated with $\lambda$. In other words, the $H$-eigenvalues of $\A$ and $P^m \A$ are the same for any permutation matrix $P\in\R^{n\times n}$. Thus (ii) is obtained. \qed

A sufficient condition is provided as below for the double nonnegativity, which will be very useful in the subsequent analysis.\vskip 2mm

\begin{proposition}\label{add} Suppose $\A$ is a nonnegative tensor and $\A=\B+\C$ with $\B=\sum\limits_{k=1}^{r_1} \left(u^{(k)}\right)^m$ and $\C=\sum\limits_{j=1}^{r_2} \left(v^{(j)}\right)^m$, where $u^{(k)}\in\R^n$ and $v^{(j)}\in\R^n_+$ for all $k\in [r_1]$ and $j\in [r_2]$. If there exists some $i_0\in [n]$ such that $u^{(k)}_{i_0}>0$ for all $k\in [r_1]$, then $\A$ is doubly nonnegative.
\end{proposition}
\noindent{\bf Proof.} The assertion is obvious for even order $m$ due to the convexity of the positive semidefinite tensor cone. Efforts are then made on the case of odd $m$. For any $H$-eigenvalue of $\A$ with its associated $H$-eigenvector $x$, we have
$$\lambda x_i^{m-1}=\left(\A x^{m-1}\right)_i=\sum\limits_{k=1}^{r_1} \left(x^{T}u^{(k)}\right)^{m-1}u_i^{(k)}+\sum\limits_{j=1}^{r_2} \left(x^{T}v^{(j)}\right)^{m-1}v_i^{(j)},~\forall i\in [n].\eqno(3.1)$$
If $x^{T}u^{(k)}=0$ for all $k\in [r_1]$, along with $x\neq 0$, we can find some $i\in [n]$ such that $x_i\neq 0$ and hence $\lambda =\frac{\sum\limits_{j=1}^{r_2} \left(x^{T}v^{(j)}\right)^{m-1}v_i^{(j)}}{x_i^{m-1}}\geq 0$. If there exists some $\bar{k}\in [r_1]$ such that $x^{T}u^{(\bar{k})}\neq 0$, (3.1) yields that $\lambda x_{i_0}^{m-1}>0$. Thus $\lambda>0$. This completes the proof.\qed

The condition as required in Proposition \ref{add} is sufficient but not necessary for the double nonnegativity as you can see in the following example.

\begin{example} Let $\A=\B+\C$ with $\B=\left((1,0,-1)^{T}\right)^3+\left((-1,0,0)^{T}\right)^3$ and $\C=\left((1,1,1)^{T}\right)^3$. It is easy to verify that $\A$ is nonnegative and $\B$ does not satisfy the condition as stated in Proposition \ref{add}. Furthermore, it is impossible to find any rank-one decomposition for $\B$ to satisfy the desired condition since $b_{111}=b_{222}=0$ and $b_{333}=-1$. However, $\A$ is still doubly nonnegative. We can show the nonnegativity of all $H$-eigenvalues of $\A$ by contrary. Assume that there exists some $H$-eigenvalue $\lambda<0$ with its associated $H$-eigenvector $x$. By definition, we have
$$\A x^2=(x_1-x_3)^2\left[
                      \begin{array}{c}
                        1 \\
                        0 \\
                        -1 \\
                      \end{array}
                    \right]+x_1^2\left[
                      \begin{array}{c}
                        -1 \\
                        0 \\
                        0 \\
                      \end{array}
                    \right]+(x_1+x_2+x_3)^2\left[
                      \begin{array}{c}
                        1 \\
                        1 \\
                        1 \\
                      \end{array}
                    \right]=\lambda \left[
                      \begin{array}{c}
                        x_1^2 \\
                         x_2^2 \\
                        x_3^2 \\
                      \end{array}
                    \right].$$ By the second equation, we have $x_2=0$ and $x_1=-x_3\neq 0$ since $x\neq 0$. Substituting these values into the first equation, it follows that $\lambda=3$ which is a contradiction to the assumption that $\lambda<0$.
\end{example}

\subsection{Double Nonnegativity of Structured Tensors}
Many even-order structured tensors have been shown to be positive semidefinite, such as the diagonally dominant tensor, the generalized diagonally dominant tensor, the $H$-tensor with nonnegative diagonal entries, the complete Hankel tensor, the strong Hankel tensor, the $MB_0$-tensor, the quasi-double $B_0$-tensor, the double $B_0$-tensor and the $B_0$-tensor. For the odd order case, we will prove that with nonnegative entries, they are all doubly nonnegative as well, as the following theorem elaborates.

\begin{theorem}\label{DNN-Tensors} Let $\A$ be a nonnegative symmetric tensor. If one of the following conditions holds

(i) $\A$ is a diagonally dominant tensor;

(ii) $\A$ is a generalized diagonally dominant tensor;

(iii) $\A$ is an $H$-tensor;

(iv) $\A$ is a complete Hankel tensor;

(v) $\A$ is a strong Hankel tensor;

(vi) $\A$ is an $MB_0$-tensor;

(vii) $\A$ is the signless Laplacian tensor of a uniform $m$-hypergraph;

\noindent then $\A$ is a doubly nonnegative tensor.

\end{theorem}

\noindent{\bf Proof.} (i) The nonnegativity of $H$-eigenvalues follows from Theorem 6 in \cite{Q2005}. (ii) By definition, we can find some positive diagonal matrix $D$ such that $\A D^{m-1}$ is diagonally dominant. This further implies that $\B:=\A D^{1-m}\underbrace{D\cdots D}\limits_{m-1}$ is diagonally dominant, and hence all $H$-eigenvalues are nonnegative. Applying Proposition \ref{invariance}, we can get (ii). (iii) The desired assertion can be obtained by invoking Theorem 4.9 in \cite{KMB2015}, together with (ii). (iv) For any complete Hankel tensor $\A$ with its Vandermonde decomposition $$\A=\sum\limits_{k=1}^{r}\alpha_k\left(u^{(k)}\right)^m,$$
where $\alpha_k>0$, $u^{(k)}=\left(1,\xi_k,\ldots, \xi_k^{n-1}\right)^{T}\in\R^n$ for all $k\in [r]$, the desired assertion follows readily from Proposition \ref{add} by setting $i_0=1$ and $\C=0$.  (v) It is known from \cite{DQW2015} that for any strong Hankel tensor $\A\in \S_{m,n}$, it has an augmented Vandermonde decomposition
$$\A=\sum\limits_{k=1}^{r-1} \alpha_k\left(u^{(k)}\right)^m+\alpha_{r}\left(e^{(n)}\right)^{m},\eqno(3.2)$$
where $\alpha_k>0$, $u^{(k)}=\left(1,\xi_k,\ldots, \xi_k^{n-1}\right)^{T}\in\R^n$, $\xi_k\in\R$, for all $k\in [r-1]$, $\alpha_r\geq 0$ and $e^{(n)}\in \R^n$ is the $n$th column of the identity matrix. Utilizing Proposition \ref{add} again, we can get the double nonnegativity of any nonnegative strong Hankel tensors. (vi) It is known from Theorem 7 in \cite{LQL2015} that for any nonnegative symmetric $MB_0$-tensor $\A$, either $\A$ is a symmetric $M$-tensor itself or we have $$\A={\mathcal{M}}+\sum\limits_{k=1}^sh_k{\mathcal{E}}^{J_k},$$ where $\mathcal{M}$ is a symmetric $M$-tensor, $s$ is a positive integer, $h_k>0$ and $J_k\subset [n]$, for $k=1,\cdots, s$. When $m$ is even, the desired result can be derived from the positive semidefiniteness. When $m$ is odd, the assertion is obvious when $\A$ is an $M$-tensor. To show for the latter case, we first claim that for any symmetric $M$-tensor ${\mathcal{M}}$ and any vector $x\in\R^n\setminus \{0\}$, there always exists some $i\in supp(x):=\{i\in [n]:x_i\neq 0\}$ such that $\left({\mathcal{M}}x^{m-1}\right)_i\geq 0$. Assume on the contrary that there exists some nonzero $x$ such that for any $i\in supp(x)$, $\left({\mathcal{M}} x^{m-1}\right)_i<0$. Let $\alpha_i=-\left({\mathcal{M}} x^{m-1}\right)_i /x^{m-1}_i$, for all $i \in supp(x)$. Obviously, $\alpha_i>0$ for all $i\in supp(x)$. Thus, $\left(\bar{\mathcal{M}}+\sum\limits_{i\in supp(x)}\alpha_i \left(e^{(i)}\right)^m\right)\bar{x}^{m-1}=0$, where $\bar{\mathcal{M}}$ is the principal subtensor of ${\mathcal{M}}$ and $\bar{x}$ the sub-vector of $x$ generated by the index set $supp(x)$. This comes to a contradiction to the fact that $\bar{\mathcal{M}}+\sum\limits_{i\in supp(x)}\alpha_i\left(e^{(i)}\right)^m$ is a strong $M$-tensor by the property of $M$-tensors. This shows our claim. Now for any $H$-eigenvalue of $\A$ with its associated $H$-eigenvector $x$, we have
$$\lambda x_i^{m-1}=\left(\A x^{m-1}\right)_i=\left({\mathcal{M}} x^{m-1}\right)_i+\sum\limits_{k=1}^s h_k\left({\mathcal{E}}^{J_k}x^{m-1}\right)_i,~\forall i\in [n]. $$
From our claim, we can find some $i\in supp(x)$ such that $\left({\mathcal{M}}x^{m-1}\right)_i\geq 0$ and hence $\lambda=\frac{\left(\A x^{m-1}\right)_i}{x^{m-1}_i} \geq 0$. (vii) By definition, the signless Laplacian tensor is a nonnegative tensor. By \cite{Qi2014}, all $H$-eigenvalues of the
signless Laplacian tensor are nonnegative. Hence, it is a doubly nonnegative tensor. This completes the proof.\qed

\begin{remark} It has been shown in \cite{LL2015,LQL2015} that $$\{B_0\text{-tensors}\}\subset \{double~B_0\text{-tensors}\}\subset \{\text{quasi-double}~ B_0\text{-tensors}\}\subset \{MB_0\text{-tensors}\}.$$ From (vii) in Theorem \ref{DNN-Tensors}, all these tensors with nonnegative entries are doubly nonnegative when they are symmetric.
\end{remark}

Similarly, we can get the following results on strongly doubly nonnegative tensors.\vskip 2mm

\begin{proposition}\label{SDNN-Tensors} Let $\A\in \S_{m,n}$ be a nonnegative symmetric tensor. If one of the following conditions holds

(i) $\A$ is a strictly diagonally dominant tensor;

(ii) $\A$ is a generalized strictly diagonally dominant tensor;

(iii) $\A$ is a strong $H$-tensor;

(iv) $\A$ is an $MB$-tensor (or $B$-tensor, or quasi-double $B$-tensor, or double $B$-tensor);

\noindent then $\A$ is a strongly doubly nonnegative tensor.
\end{proposition}



\begin{proposition}\label{circulant2} Let $\A\in\S_{m,n}$ be a nonnegative circulant tensor with its root tensor $\A_1$ and its diagonal entry $c_0$. If $2c_0-\A_1 e^{m-1}\geq 0$, then $\A$ is doubly nonnegative.
\end{proposition}
\noindent{\bf Proof.} By Definition \ref{def-circulant}, the condition $2c_0-\A_1 e^{m-1}\geq 0$ implies that $\A$ is diagonally dominant. Thus, Theorem \ref{DNN-Tensors} immediately leads to the desired result. \qed

For circulant tensors, we also get the following properties on copositivity, which can be regarded as necessary conditions for the double nonnegativity.\vskip 2mm

\begin{theorem}\label{circulant1} Let $\A\in\T_{m,n}$ be a circulant tensor with its root tensor $\A_1$. If $\A_1$ is copositive (strictly copositive, respectively), then $\A$ is copositive (strictly copositive, respectively). Moreover, if $\A$ is a doubly circulant tensor, then $\A$ is copositive (strictly copositive, respectively) if and only if $\A_1$ is copositive (strictly copositive, respectively).

\end{theorem}
\noindent{\bf Proof.} Let $a_{i_1\ldots i_m}$ be the $(i_1,\ldots, i_m)$th entry of $\A$ and $A_k$ ($k\in[n]$) be its row tensors. Invoking Proposition 2 in \cite{CQi2015}, we have $\A_{k+1}=P^m(\A_k)$, where $P=(p_{ij})\in\R^{n\times n}$ is a permutation matrix with $p_{i+1 i}=1$ for $i\in[n-1]$, $p_{1n}=1$ and $p_{ij}=0$ otherwise. If $\A_1$ is copositive, then for any $x\in\R^n_+$, $\underbrace{P\cdots P}\limits_{k-1}x\in\R^n_+, \forall k\in [n]$ and
$$  \A x^m = \sum\limits_{k=1}^n x_k\A_kx^{m-1}=\sum\limits_{k=1}^n x_k \A_1 \big(\underbrace{P\cdots P}\limits_{k-1}\big)^{m-1}x^{m-1}= \sum\limits_{k=1}^n x_k \A_1 \big(\underbrace{P\cdots P}\limits_{k-1}x\big)^{m-1} \geq 0,$$
where the last inequality follows from the copositivity of $\A_1$. Since $\underbrace{P\cdots P}\limits_{k-1}x\in\R^n_+\setminus\{0\}, \forall k\in [n]$ for any $x\in\R^n_+\setminus \{0\}$, we can similarly show the case of strict copositivity of $\A$ if $\A_1$ is strictly copositive. To get the moreover part, it suffices to show the necessity. By the definition of doubly circulant tensor, we have all row tensors $\A_k$ coincide with its root tensor $\A_1$. Thus
$$\A x^{m}=\sum\limits_{k=1}^nx_k\A_k x^{m-1}=\left(\sum\limits_{k=1}^n x_k\right) \A_1 x^{m-1}, \forall x\in \R^n.$$
If $\A$ is copositive, for any $x\in\R^n_+\setminus \{0\}$, $\sum\limits_{k=1}^n x_k >0$. Thus, $\A_1 x^{m-1}= \left(\A x^m\right)/\left(\sum\limits_{k=1}^n x_k\right)\geq 0$. This shows the copositivity of $\A_1$. Similarly, we can prove the strictly copositive case. \qed
\vskip 2mm

\section{Completely Positive Tensors}
Completely positive tensors form an important subclass of doubly nonnegative tensors. In this section, many useful properties are explored for this special class of tensors. We start by the dominance properties which serve as a powerful tool for excluding many nonnegative tensors from the set of completely positive tensors.\vskip 2mm

\subsection{Dominance Properties}

Dominance properties of completely positive tensors were studied in \cite{QXX2014}.

\begin{definition} A tensor $\A\in \S_{m,n}$ is said to have the \emph{zero-entry dominance property} if $a_{i_1\ldots i_m}=0$ implies that $a_{j_1\ldots j_m}=0$ for any $(j_1,\ldots, j_m)$ satisfying $\{j_1,\ldots, j_m\}\supseteq\{i_1,\ldots, i_m\}$.
\end{definition}

\begin{proposition}\label{zero} If $\A$ is a completely positive tensor, then $\A$ has the zero-entry dominance property.
\end{proposition}
\noindent{\bf Proof.} This follows directly from Theorem 3 in \cite{QXX2014}.\qed



Utilizing the zero-entry dominance property, we can exclude some doubly nonnegative tensors from the class of complete positive tensor very efficiently.

\begin{proposition}\label{signless} The signless Laplacian tensor of a nonempty uniform $m$-hypergraph for $m \ge 3$  is not completely positive.
\end{proposition}
\noindent{\bf Proof.} Suppose that $m \ge 3$ and $G$ is a nonempty uniform $m$-hypergraph. Suppose that $(j_1, \cdots, j_m)$ is an edge of $G$.  Let
$\Q = (q_{i_1\cdots i_m})\in \S_{m,n}$ be the signless Laplacian tensor of $G$.  By definition, $q_{j_1\ldots j_m} = {1 \over (m-1)!} \neq 0$. Note that $q_{j_1j_1\ldots j_1j_2}=0$ by the definition of signless Laplacian tensors. Obviously, the zero-entry dominance property fails and hence $\Q$ is not completely positive. \qed


\begin{proposition}\label{Hankel0} Let $\A\in \S_{m,n}$ be a Hankel tensor, and $v=(v_0,\ldots, v_{(n-1)m})^T\in \R^{(n-1)m+1}$ be its generating vector.

(i) If $v_{0}=v_{(n-1)m}=0$, then $\A\in CP_{m,n}$ if and only if $\A=O$;

(ii) If $\A\in CP_{m,n}$ and $v_{(i-1)m}=0$ for some $2\leq i\leq n-1$, then $v_{0}\geq 0$, $v_{(n-1)m}\geq 0$ and $\A=v_0 e_1^{m}+v_{(n-1)m}e_n^{m}$;

(iii) If $v_0=0$ and $v_j\neq 0$ for some $j\in [m-1]$, then $\A$ is not completely positive.

\end{proposition}
\noindent{\bf Proof.} (i) It is trivial that $O$ is a completely positive tensor. If $\A$ is completely positive and $v_0=v_{(n-1)m}=0$, then $a_{1i_2\ldots i_m}=0$ for all $i_2$, $\ldots$, $i_m\in [n]$. Thus $v_j=0$, for any $j\in [(m-1)n+1-m]$. When $n$=2, then $\A=0$. When $n\geq 3$, then $n\geq 2+\frac{1}{m-1}$, which implies that $m\geq (m-1)n+1-m$. Thus $a_{2\ldots 2}=0$. Using the zero-entry dominance property again, we get $v_j=0$ for all $j\in [(m-1)n+2-m]$. Keep on doing this, we can find that for any given $k\in [n-1]$, $v_j=0$ for all $j\in [(m-1)n+k-1-m]$, then $a_{k\ldots k}=0$ since $n\geq k+\frac{1}{m-1}$. The zero-entry dominance property and the fact $v_{(n-1)m}=0$ finally give us $\A=0$. Thus (i) is obtained. Using the similar proof as in (i), we can prove that $a_{k\ldots k}=0$ for all $k=2,\ldots,n-1$. Thus, the zero-entry dominance property shows $\A=v_0 e_1^{m}+v_{(n-1)m}e_n^{m}$. By the nonnegativity of $\A$, $v_0$ and $v_{(n-1)m}$ are nonnegative. This implies the assertion in (ii). By definition, we know that $a_{1\ldots 1}=0$. Theorem \ref{zero} tells us that $a_{i_1\ldots i_m}=0$ for any $(i_1,\ldots, i_m)$ satisfying $1\in \{i_1,\ldots, i_m\}$. Note that if $i_1+\cdots+i_m\leq 2m-1$, then $1 \in \{i_1,\ldots, i_m\}$. This leads to (iii).  \qed

\begin{proposition}\label{circulant0} Let $\A$ be a Toeplitz tensor with its diagonal entry $0$. Then $\A$ is completely positive if and only if $\A=O$.
\end{proposition}
\noindent{\bf Proof.} The sufficiency is trivial. If $\A$ is completely positive and $a_{1\ldots 1}=0$, by the definition of Toeplitz tensors, we have $a_{i\ldots i}=0$ for all $i\in [n]$. Invoking the zero-entry dominance property in Theorem \ref{zero}, it follows that $\A=O$. \qed

The zero-entry dominance property may work very well for excluding some doubly nonnegative tensors with zero entries from the class of completely positive tensors. The strong dominance property, as described below, may work for the case of some other doubly nonnegative tensors that may have all positive entries.

Let $I=\left\{\left(i_1^{(1)},\ldots, i_m^{(1)}\right),\ldots, \left(i_1^{(s)},\ldots, i_m^{(s)}\right)\right\}$ with $\left\{i_1^{(p)},\ldots, i_m^{(p)}\right\}\subseteq \{j_1,\ldots, j_m\}$ for any $p\in [s]$. For any given index $i\in \{j_1,\ldots, j_m\}$, if it appears $t$ times in $\{j_1,\ldots, j_m\}$, then it appears in $I$ $st$ times. Then we call $I$ an $s$-duplicate of $\left(j_1,\ldots, j_m\right)$. 


\begin{proposition}[Strong Dominance,~Theorem 4, \cite{QXX2014}]\label{strong dominance property} Suppose that $\A=\left(a_{i_1\ldots i_m}\right)\in\S_{m,n}$ is completely positive and $I=\left\{\left(i_1^{(1)},\ldots, i_m^{(1)}\right),\ldots, \left(i_1^{(s)},\ldots, i_m^{(s)}\right)\right\}$ is an $s$-duplicate of $\left(j_1,\ldots, j_m\right)$ for some given $j_l\in [n]$, $l\in [m]$. Then $s^{-1}\sum\limits_{p=1}^s a_{i_1^{(p)}\cdots i_m^{(p)}}\geq a_{j_1\ldots j_m}$.
\end{proposition}

The aforementioned strong dominance property provides us a way to exclude some positive doubly nonnegative tensors from the class of completely positive tensors, as you will see in Section 7.

\subsection{Spectral Properties}

It is known from \cite{QXX2014} that completely positive tensors have the following spectral properties.

\begin{theorem}[Theorems 1,2 \cite{QXX2014}] \label{th1} Let $\A\in\S_{m,n}$ be a completely positive tensor. Then

(i) all $H$-eigenvalues of $\A$ are nonnegative;

(ii) When $m$ is even, all its $Z$-eigenvalues are nonnegative; When $m$ is odd, a $Z$-eigenvector associated with a positive (negative) $Z$-eigenvalue of $\A$ is nonnegative (nonpositive).

\end{theorem}

For strongly completely positive tensors, we have the following spectral properties.

\begin{theorem}\label{th2} Let $\A\in\S_{m,n}$ be a strongly completely positive tensor. Then

(i) all $H$-eigenvalues of $\A$ are positive;

(ii) all its $Z$-eigenvalues are nonzero. Moreover, when $m$ is even, all its $Z$-eigenvalues are positive, when $m$ is odd, a $Z$-eigenvector associated with a positive (negative) $Z$-eigenvalue of $\A$ is nonnegative (nonpositive).
\end{theorem}

\noindent{\bf Proof.} Write $\A$ as $\A=\sum\limits_{k=1}^r \left(u^{(k)}\right)^m$, where $u^{(k)}\in \R^n_+$ and $$span\{u^{(1)},\ldots, u^{(r)}\}=\R^n.\eqno(4.1)$$  (i) Assume on the contrary that $\A$ has $\lambda =0$ as one of its $H$-eigenvalues, and the corresponding $H$-eigenvector is $x$. Certainly $x\neq 0$. When $m$ is even, by the definition of $H$-eigenvalue, we have
$$ 0=\lambda\sum\limits_{i=1}^m x_i^m=\A x^{m}=\sum\limits_{k=1}^r\left(x^{T}u^{(k)}\right)^m.$$
The nonnegativity of each term in the summation on the right hand side immediately leads to $x^{T}u^{(k)}=0$ for all $k\in [r]$. Invoking the condition in (4.1), $x$ has no choice but $0$, which comes to a contradiction since $x$ is an $H$-eigenvector. Thus, all $H$-eigenvalues of $\A$ is positive when the order is even. When $m$ is odd, it is known by definition that
$$0=\lambda x_i^{m-1}=\left(\A x^{m-1}\right)_i=\sum\limits_{k=1}^r\left(x^{T}u^{(k)}\right)^{m-1}u^{(k)}_i, ~\forall i\in[n]. \eqno(4.2)$$ Together with the involved nonnegativity of each term in the summation on the right hand side, (4.2) implies that
$$\left(x^{T}u^{(k)}\right)^{m-1}u^{(k)}_i=0,  ~\forall i\in[n].\eqno(4.3)$$ In addition, the condition (4.1) implies that we can pick $n$ vectors from the set $\{u^{(1)},\ldots, u^{(r)}\}$ to span the whole space $\R^n$. Without loss of generality, let's say they are $u^{(1)}$, $\ldots$, $u^{(n)}$. Trivially, for any $k\in [n]$, $u^{(k)}\neq 0$. Therefore, there always exists an index $i_k\in [n]$ such that $u^{(k)}_{i_k}\neq 0$. Thus (4.3) implies that $x^{T}u^{(k)}=0$, for all $k\in [n]$. This immediately leads to $x=0$. The same contradiction arrives and hence all $H$-eigenvalues of $\A$ should be positive when $m$ is odd.
(ii) Assume on the contrary that $\lambda =0$ is a $Z$-eigenvalue with $x$ as its $Z$-eigenvector. When $m$ is even, it follows by definition that
$$0=\lambda x^{T}x=\A x^m=\sum\limits_{k=1}^k \left(x^{T}u^{(k)}\right)^m.$$
Since each term $\left(x^{T}u^{(k)}\right)^m\geq 0$ ($k\in [r]$), it follows readily that $x^{T}u^{(k)}=0$ for all $k\in [r]$. In addition, together with (4.1), we get $x=0$, which contradicts to the assumption that $x$ is a $Z$-eigenvector. When $m$ is odd, it comes directly that
$$0=\lambda x_i=\A x^m=\sum\limits_{k=1}^r\left(x^{T}u^{(k)}\right)^{m-1} u^{(k)}_i, ~\forall i\in [n].$$
The nonnegativity of each term $\left(x^{T}u^{(k)}\right)^{m-1} u^{(k)}_i$ for all $k\in [r]$ and all $i\in [n]$ yields
$$\left(x^{T}u^{(k)}\right)^{m-1} u^{(k)}_i=0, \forall k\in [r], ~\forall i\in [n].\eqno(4.4)$$
Pick up $n$ linearly independent vectors from the set $\{u^{(1)},\ldots, u^{(r)}\}$, simply say $u^{(1)},\ldots, u^{(n)}$. By the observation that $u^{(k)}\neq 0$ for all $k\in [n]$, we can always find some index $i_k\in [n]$ such that $u^{(k)}_{i_k}\neq 0$, for all $k\in [n]$. Substituting this into (4.4), we can obtain $x^{T}u^{(k)}=0$ for all $k\in [n]$. This gives us $x=0$, which is a contradiction since $x$ is a $Z$-eigenvector. This completes the proof. \qed

\subsection{Necessary and/or Sufficient Conditions}
Some necessary and/or sufficient conditions for (strongly) completely positive tensors are presented in this subsection.

\begin{proposition}\label{cp1} For any given nonnegative matrix $P\in\R^{l\times n}$, if $\A$ is completely positive, then $P^m(\A)$ is also completely positive.
\end{proposition}
\noindent{\bf Proof.} Let $\A=\sum\limits_{k=1}^r \left(u^{(k)}\right)^m$, where $u^{(k)}\in \R^n_+$. By employing the decomposition invariance in Theorem 2.2 of \cite{LQY2015}, we know that $P^m(\A)=\sum\limits_{k=1}^r \left(Pu^{(k)}\right)^m$. If $P$ is nonnegative, then so are $Pu^{(k)}$ for all $k\in [r]$. \qed

\begin{proposition}\label{nonsingular} For any nonnegative nonsingular $P\in\R^{n\times n}$, $\A$ is (strongly) completely positive  if and only if $P^m(\A)$ is (strongly) completely positive.
\end{proposition}
\noindent{\bf Proof.} Let $\A=\sum\limits_{k=1}^r \left(u^{(k)}\right)^m$, where $u^{(k)}\in \R^n_+$. By employing the decomposition invariance in Theorem 2.2 of \cite{LQY2015}, we know that $P^m(\A)=\sum\limits_{k=1}^r \left(Pu^{(k)}\right)^m$. The nonsingularity of $P$ implies that $span\{Pu^{(1)},\ldots, Pu^{(r)}\}= span\{u^{(1)},\ldots, u^{(r)}\}$, and the nonnegativity of $P$ implies that $Pu^{(k)}\in \R^n_+$ for all $k\in [r]$. Conversely if $\B:=P^m(\A)$ is completely positive (definite), then similarly we can prove that $\A=\left(P^{-1}\right)^m \B$ is completely positive (definite).\qed

\begin{proposition}\label{cp3} Let $\A\in\S_{m,n}$ be a completely positive tensor. Then

(i) for any even integer $l\in [m]$, $\A x^l\in \S_{m-l,n}$ is also completely positive for any $x\in\R^n$;

(ii) for any integer $t\in [m]$, $\A x^t\in \S_{m-t,n}$ is also completely positive for any $x\in\R^n_+$;

(iii) for any $\Gamma\subseteq [n]$, the principal subtensor $\A_{\Gamma}\in \S_{m, |\Gamma|}$ is also completely positive.
\end{proposition}
\noindent{\bf Proof.} Let $\A=\sum\limits_{k=1}^r \left(u^{(k)}\right)^m$, where $u^{(k)}\in \R^n_+$. It follows that for any $x\in\R^n$ and any integer $l\in [m]$,
$$\A x^l=\sum\limits_{k=1}^r \left(x^{T}u^{(k)}\right)^l \left(u^{(k)}\right)^{m-l}.$$ Thus (i) and (ii) can be obtained by definition. (iii) is from Proposition 2 in \cite{QXX2014}. \qed

Similarly, we can get the following properties for strongly completely positive tensors.

\begin{proposition}\label{cp4} Let $\A\in\S_{m,n}$ be a strongly completely positive tensor. Then
\begin{itemize}
\item[(i)] for any integer $l\in [m]$, $\A x^l\in \S_{m-l,n}$ is also strongly completely positive for any $x\in\R^n_{++}$;

\item[(ii)] for any $\Gamma\subseteq [n]$, the principal subtensor $\A_{\Gamma}\in \S_{m, |\Gamma|}$ is also strongly completely positive.
\end{itemize}
\end{proposition}

The Hadamard product preserves the complete positivity as shown in Proposition 1 in \cite{QXX2014}. It also preserves the strong complete positivity as stated below.

\begin{proposition}\label{cp5} Let $\A$, $\B\in\S_{m,n}$. If $\A$ and $\B$ are strongly completely positive, then $\A\circ \B$ is strongly completely positive.
\end{proposition}
\noindent{\bf Proof.} We first claim that if $U=\left(u^{(1)}~u^{(2)}~\ldots ~u^{(n)}\right)$, $V=\left(v^{(1)}~v^{(2)}~\ldots ~v^{(n)}\right)$ are any two nonsingular matrices in $\R^{n\times n}$, then $$span\{u^{(1)}\circ v^{(1)},u^{(1)}\circ v^{(2)},\ldots, u^{(1)}\circ v^{(n)},u^{(2)}\circ v^{(1)}, u^{(2)}\circ v^{(2)}, \ldots, u^{(n)}\circ u^{(n)} \}=\R^n.\eqno(4.5)$$ The nonsingularity of $U$ indicates that $u^{(1)},u^{(2)},\ldots,u^{(n)}$ can form a basis for $\R^n$. Thus, we can find $a_{ik}$, $i$, $k\in [n]$, such that
$$e^{(i)}=\sum\limits_{k=1}^n a_{ik}u^{(k)},~\forall i\in [n],\eqno(4.6)$$
where there exists at least one nonzero element among $a_{i1}$, $\ldots$, $a_{in}$ for any $i\in [n]$. The equalities in (4.6) derive that
$$e^{(i)} \circ v^{(j)}=\sum\limits_{k=1}^n a_{ik}\left(u^{(k)}\circ v^{(j)}\right),~\forall i,~j\in[n]. \eqno(4.7)$$
By the nonsingularity of $V$, it follows that
 $0\neq \det(V)=\sum_{\sigma\in S_n} sgn(\sigma)\Pi_{i=1}^n v^{(i)}_{\sigma_i}$, where $\sigma=(\sigma_1,\ldots, \sigma_n)$ is a permutation of $[n]$, $ sgn(\sigma)$ is the signature of $\sigma$, 
 $S_n$ is the set of all permutations of $[n]$. Thus, there always exists some permutation $\sigma$ such that $\Pi_{i=1}^n v^{(i)}_{\sigma_i}\neq 0$. Without loss of generality, we assume that $\sigma=(1,2,\ldots, n)$, that is, $v^{(i)}_i\neq 0$ for all $i\in [n]$. This together with (4.7) yields that
 $$v^{(i)}_ie^{(i)}=\sum\limits_{k=1}^n a_{ik}\left(u^{(k)}\circ v^{(i)}\right), ~\forall i\in[n],$$
 which indicates that $e^{(i)}=\sum\limits_{k=1}^n \frac{a_{ik}}{v^{(i)}_i}\left(u^{(k)}\circ v^{(i)}\right)$, for all $i\in[n]$. Thus, $\{e^{(1)},\ldots, e^{(n)}\}$ can be linearly expressed by $\{u^{(i)}\circ v^{(j)}\}_{i,j=1,\ldots, n}$, and our claim is proven. Now we consider any two strongly completely positive tensors $\A$ and $\B$ with their corresponding nonnegative rank-one decompositions $\A=\sum\limits_{i=1}^r \left(u^{(i)}\right)^m$, and $\B=\sum\limits_{j=1}^{r^{'}} \left(v^{(j)}\right)^m,$ where $span\{u^{(1)},\ldots, u^{(r)}\}=span\{v^{(1)},\ldots, v^{(r^{'})}\}=\R^n$. Easily we can verify that $\A\circ \B=\sum\limits_{i=1}^r\sum\limits_{j=1}^{r^{'}}\left(u^{(i)}\circ v^{(j)}\right)^m$. The involved $u^{(i)}\circ v^{(j)}$ is certainly nonnegative by the nonnegativity of $u^{(i)}$ and $v^{(j)}$ for all $i$, $j\in [n]$. Note that $r$ and $r^{'}$ should be no less than $n$. Therefore, we can always pick up $n$ vectors from $\{u^{(1)},\ldots, u^{(r)}\}$ to form a basis of $\R^n$. Let's simply say they are $u^{(1)},\ldots, u^{(n)}$. Similarly, we can do this to $v^{(1)},\ldots, v^{(r^{'})}$ and get $n$ linearly independent vectors, namely $v^{(1)},\ldots, v^{(n)}$. The aforementioned claim tells us that all involved vectors  $u^{(i)}\circ v^{(j)}$, $i$, $j\in [n]$, can span the whole space $\R^n$, which means $\A\circ \B$ is strongly completely positive.\qed



\section{Cones of Doubly Nonnegative Tensors}
In this section, many tensor cones and their relationship in the class of doubly nonnegative tensors are discussed. Particularly, the gap existing between the doubly nonnegative tensor cone and the completely positive tensor cone is analyzed. Following from the matrix case \cite{BM2003}, the set of all doubly nonnegative tensors of order $m$ and dimension $n$ is denoted by $DNN_{m,n}$. Similarly, we use $SDNN_{m,n}$ to denote the set of all strongly doubly nonnegative tensors in $DNN_{m,n}$. Adopting the notations in the literature (see e.g., \cite{QXX2014}), $CP_{m,n}$ and $COP_{m,n}$ are used to denote the sets of all completely positive tensors and all copositive tensors of order $m$ and dimension $n$, respectively. In addition, $SCP_{m,n}$ and $SCOP_{m,n}$ are used to stand for the set of all strongly completely positive tensors and the strictly copositive tensors, respectively. The set of all symmetric positive semidefinite (positive definite) tensors is denoted by $PSD_{m,n}$ ($PD_{m,n}$) for convenience.

\begin{proposition}\label{dnn1} If $m$ is even, then $DNN_{m,n}=PSD_{m,n}\cap N_{m,n}$, $SDNN_{m,n}=PD_{m,n}\cap N_{m,n}$, $int(DNN_{m,n})=PD_{m,n}\cap N^{+}_{m,n}\subseteq SDNN_{m,n}$, where $N^{+}_{m,n}:=\{\A\in\T_{m,n}: a_{i_1\ldots i_m}>0, ~i_1,\ldots, i_n\in [n]\}$.
\end{proposition}
\noindent{\bf Proof.} By definition, the desired assertions follow from Theorem 5 in \cite{Q2005}.\qed

\begin{proposition}\label{cone1} $CP_{m,n}$ and $COP_{m,n}$ are closed convex cones and they are dual to each other. Moreover,

(i)  $CP_{m,n}\subseteq DNN_{m,n}\subseteq COP_{m,n}$;

(ii) $SCP_{m,n}\subseteq SDNN_{m,n} \subseteq SCOP_{m,n}$.

\noindent Furthermore, let $\Gamma \subseteq [n]$ and $I_{\Gamma}\in\R^{n\times n}$ be the matrix with the $(i,i)$th entry $1$ if $i\in\Gamma$ and $0$ elsewhere. Then

(iii) $\left(I_{\Gamma}\right)^m CP_{m,n} \unlhd CP_{m,n}$;

(iv) $\left(I_{\Gamma}\right)^m DNN_{m,n} \unlhd DNN_{m,n}$, when $m$ is even.

\end{proposition}

\noindent{\bf Proof.} The first part is from Theorem 5 in \cite{QXX2014}. (i) can be easily verified by (i) in Theorem \ref{th1}. For (ii), the inclusion $SCP_{m,n}\subseteq SDNN_{m,n}$ can be derived from (ii) in Theorem \ref{th2}. To get the remaining inclusion $SDNN_{m,n}\subseteq SCOP_{m,n}$, suppose $\A$ is a strongly
doubly nonnegative tensor. By definition, all $H$-eigenvalues of $\A$ are positive, which implies that there exists no no-positive $H^+$-eigenvalues. This further shows that all principal subtensors of $\A$ have no no-positive $H^{++}$-eigenvalues. Applying Theorem 4.2 in \cite{SongQi2015}, $\A$ is strictly copositive.
To get (iii), we first claim that $\left(I_{\Gamma}\right)^m CP_{m,n} \subseteq CP_{m,n}$. This follows from the fact that every principal subtensor of a completely positive tensor is completely positive (Proposition 2, \cite{QXX2014}). Evidently, for any $\C \in \left(I_{\Gamma}\right)^m CP_{m,n}$, we can always find some $w^{(1)}$, $\ldots$, $w^{(l)}\in\R^n_+$ such that
$$ \A=\sum\limits_{i=1}^l \left(w^{(i)}\right)^m,~~\left(w^{(i)}\right)_j=0, ~\forall j\notin \Gamma. \eqno(5.1)$$
Now it remains to show that the following implication holds:
$$\A,~\B\in CP_{m,n},~\A+\B\in \left(I_{\Gamma}\right)^m CP_{m,n}\Rightarrow \A,~\B\in \left(I_{\Gamma}\right)^m CP_{m,n}.$$
Write
$$\A=\sum\limits_{k=1}^r \left(u^{(k)}\right)^m,~ \B=\sum\limits_{j=1}^{r^{'}} \left(v^{(j)}\right)^m, ~u^{(k)},~v^{(j)}\in\R^n_+, \forall k\in[r], j\in [r^{'}].$$
If $\A+\B\in  \left(I_{\Gamma}\right)^m CP_{m,n}$, invoking (5.1), it follows readily that
$$\sum\limits_{k=1}^ru^{(k)}_{i_1}\cdots u^{(k)}_{i_m}+\sum\limits_{j=1}^{r^{'}}v^{(j)}_{i_1}\cdots v^{(j)}_{i_m} =0, ~\forall i_1,\ldots, i_m\notin \Gamma.\eqno(5.2)$$
By setting $i_1=i_2=\ldots=i_m=t$ for any $t\notin \Gamma$, (5.2) implies that $u^{(k)}_t=v^{(j)}_t=0$ for all $k\in [r]$, $j\in[r^{'}]$ and $t\notin \Gamma$. Invoking the property of $\left(I_{\Gamma}\right)^m CP_{m,n}$ as shown in (5.1), we immediately get $\A$, $\B\in \left(I_{\Gamma}\right)^m CP_{m,n}$. Thus $\left(I_{\Gamma}\right)^m CP_{m,n}$ is a face of $CP_{m,n}$ for any $\Gamma\subseteq [n]$. Similarly, we can prove (iv). This completes the proof. \qed
\vskip 2mm
It is known that all strongly completely positive tensors with even $m$ are positive definite tensors by invoking Theorem \ref{th1} together with Theorem 5 in \cite{Q2005}. A natural question arises: Is it correct that any tensor with complete positivity and positive definiteness should be strongly completely positive? The following proposition answers it in an affirmative way.

\begin{proposition}\label{ppm} Suppose $m$ is even. Then $CP_{m,n}\cap PD_{m,n}=SCP_{m,n}$.
\end{proposition}
\noindent{\bf Proof.} First we claim that $CP_{m,n}\cap PD_{m,n}\subseteq SCP_{m,n}$. For any $\A\in CP_{m,n}\cap PD_{m,n}$, write its nonnegative rank-one decomposition as $\A=\sum\limits_{k=1}^r \left(u^{(k)}\right)^m$ with $u^{(k)}\in \R^n_+$, for all $k\in [r]$. Assume on the contrary that $\A\notin SCP_{m,n}$, which means $span\{u^{(1)},\ldots, u^{(r)}\}\neq \R^n$. Thus, there exists an $x\in\R^n\setminus \{0\}$ such that $x^{T}u^{(k)}=0$ for all $k\in [r]$. This immediately gives us $\A x^m=\sum\limits_{k=1}^r \left(x^{T}u^{(k)}\right)^m=0$, which is actually a contradiction to the hypothesis that $\A\in PD_{m,n}$. Thus the claim is proven. On the other hand, for any $\A\in SCP_{m,n}$. Apparently $\A\in CP_{m,n}$ and all $H$-eigenvalues of $\A$ are positive by applying (ii) of Theorem \ref{th1}. This positivity leads to $\A\in PD_{m,n}$ by Theorem 5 in \cite{Q2005}. Thus $A\in CP_{m,n}\cap PD_{m,n}$ and hence $SCP_{m,n}\subseteq CP_{m,n}\cap PD_{m,n}$ by the arbitrariness of $\A$. This completes the proof. \qed

Note that $PD_{m,n}$ is empty in the odd order case. However, inspired by Proposition \ref{dnn1}, it yields that $CP_{m,n}\cap PD_{m,n}=CP_{m,n}\cap SDNN_{m,n}$ when $m$ is even. Based on this observation, we can generalize the result in the above proposition to any odd order case, as the following theorem shows.

\begin{theorem}\label{SDNN-SCP} $CP_{m,n}\cap SDNN_{m,n}=SCP_{m,n}$.
\end{theorem}
\noindent{\bf Proof.} Invoking Propositions \ref{dnn1} and \ref{ppm}, we only need to consider the odd order case. For any given $\A\in CP_{m,n}\cap SDNN_{m,n}$ with its nonnegative rank-one decomposition $\A=\sum\limits_{k=1}^r \left(u^{(k)}\right)^m$. Assume on the contrary that $\A\notin SCP_{m,n}$, which means $span\{u^{(1)},\ldots, u^{(r)}\}\neq \R^n$. It is immediate to get some nonzero $x\in\R^n$ such that $x^{T}u^{(k)}=0$, for all $k\in [r]$. This shows that
$$ \left(\A x^{m-1}\right)_i=\sum\limits_{k=1}^r\left(x^{T}u^{(k)}\right)^{m-1}u_i^{(k)}=0,~\forall i\in[n],$$
which indicates that $0$ is an $H$-eigenvalue of $\A$. This contradicts to the hypothesis that $\A\in SDNN_{m,n}$. Thus, $CP_{m,n}\cap SDNN_{m,n}\subseteq SCP_{m,n}$. The reverse inclusion follows directly from (ii) of Theorem \ref{th1}. This completes the proof. \qed



It has been shown in \cite{Q2013} that the interior of the copositive tensor cone is exactly the strictly copositive tensor cone, and the interior of $DNN_{m,n}$ is the cone consisting of all strongly doubly nonnegative tensors with positive entries. For the completely positive tensor cone, when its reduces to the matrix case, i.e., $m=2$, the interior has been well characterized in \cite{Dickinson2010,DS2008}. Then how about the interior of the completely positive tensor cone in general higher order case? Some properties about the interior are discussed as follows.

\begin{proposition}\label{interior}  $int(CP_{m,n})\subseteq SCP_{m,n}\cap N^{+}_{m,n}$.
\end{proposition}
\noindent{\bf Proof.} 
The inclusion $int(CP_{m,n})\subseteq N^+_{m,n}$ is obvious. For the remaining part, the even order and the odd order cases will be analyzed separately. For the even order $m$, the inclusion $CP_{m,n}\subseteq DNN_{m,n}$ as shown in (i) of Proposition \ref{cone1}, together with Proposition \ref{dnn1}, yields that $int(CP_{m,n})\subseteq PD_{m,n}$. Invoking Proposition \ref{ppm}, we have
$$int(CP_{m,n})=int(CP_{m,n})\cap CP_{m,n} \subseteq PD_{m,n}\cap CP_{m,n}=SCP_{m,n}.$$
For the odd order $m$, let $\A=\sum\limits_{k=1}^r \left(u^{(k)}\right)^m$ be any tensor in $int(CP_{m,n})$. It is known by the definition of interior that there exists some scalar $\epsilon>0$ such that
$$\A+\epsilon \B\in CP_{m,n},~\forall \B\in {\mathbb{B}}_1:=\{\B\in \S^{m,n}: \langle \B,\B\rangle\leq 1\}.$$
Thus for any nonzero $x\in\R^n$,
$$\left((\A+\epsilon \B)x^{m-1}\right)_i=\left(\A x^{m-1}\right)_i+\sum\limits_{k=1}^l \left(x^{T}v^{(k)}_{\B}\right)^{m-1}[v^{(k)}_{\B}]_i\geq 0, \forall i\in [n],$$ where $\B=\sum\limits_{k=1}^l \left(v^{(k)}_{\B}\right)^{m}$ with all $v^{(k)}_{\B}\in\R^n$ for any given $\B\in {\mathbb{B}}_1$. By the arbitrariness of $\B$, it is evident that for any nonzero $x$, there exists some $j\in [n]$ such that $\left(\A x^{m-1}\right)_j>0$ and $x_j\neq 0$. In this regard, for any $H$-eigenpair $(\lambda,x)$ of $\A$, due to the definition, we have $\left(\A x^{m-1}\right)_i=\lambda_i x_i^{m-1}$, for all $i\in [n]$. Combining with the aforementioned argument, there exists some $j\in [n]$ such that $\left(\A x^{m-1}\right)_j>0$ and $x_j\neq 0$. This leads to $\lambda>0$. Thus $\A\in SDNN_{m,n}$ and hence $int(CP_{m,n})\subseteq SDNN_{m,n}$. Utilizing Theorem \ref{SDNN-SCP}, it follows that
$$int(CP_{m,n})=int(CP_{m,n})\cap CP_{m,n}\subseteq SDNN_{m,n} \cap CP_{m,n}=SCP_{m,n}.$$
  \qed



It is worth pointing out that $int(CP_{m,n})$ is a proper subset of $SCP_{m,n}\cap N^{+}_{m,n}$ with $m\geq 3$, as the following example shows.

\begin{example}\label{bad} Let $\A:=\I+e^3$. It is easy to verify that $\A\in SCP_{3,3}\cap N^+_{3,3}$. Let $\B\in\S_{3,3}$ with $b_{113}=b_{131}=b_{311}=b_{223}=b_{232}=b_{322}=1$, $b_{123}=b_{132}=b_{213}=b_{231}=b_{312}=b_{321}=-1$, and others zero. Then $\B x^3=3x_3(x_1-x_2)^2\geq 0$ for all $x=(x_1,x_2,x_3)^T\in\R^3_+$. Thus $\B\in COP_{3,3}$. However, $\langle \A,\B\rangle=0$. Thus $\A\notin int(CP_{3,3})$.
\end{example}

In \cite{Dickinson2010}, Dickinson has proposed an improved characterization of the interior of the completely positive cone for matrices, i.e., $m=2$. Following his work, the interior can be characterized as $$int(CP_{2,n})=\left\{A\in \S_{2,n}: \begin{subarray}{l}
A=\sum\limits_{k=1}^r u^{k}(u^{k})^{T},u^{(k)}\in\R^n_+, \forall i\in[r],\\
 u^{(1)}\in\R^n_{++}, span\{u^{(1)},\ldots, u^{(r)}\}=\R^n
 \end{subarray} \right\}.$$
A natural question is: Can we generalize this characterization to higher-order tensors and use the set
 $$INT=\left\{\A\in \S_{m,n}: \begin{subarray}{l}
\A=\sum\limits_{k=1}^r (u^{k})^{m},u^{(k)}\in\R^n_+, \forall i\in[r],\\
 u^{(1)}\in\R^n_{++}, span\{u^{(1)},\ldots, u^{(r)}\}=\R^n
 \end{subarray} \right\} $$
to describe $int(CP_{m,n})$ for general $m$ and $n$?  Unfortunately, the answer is no, as you can see the tensor $\A$ defined in Example \ref{bad}, where $\A \in INT\setminus int(CP_{m,n})$. In \cite{PVZ2014}, the sum of rank-one tensors, served as the standard basis of $\S_{m.n}$, is pointed out to be a special interior point in $CP_{m,n}$. But the full characterization of the entire interior is unknown.

\noindent{\bf Question 1:} How to characterize $int(CP_{m,n})$ for general $m$ and $n$?


In \cite{Dickinson2010}, it is shown that for any copositive matrix $A$, if there exists some positive $x\in\R^n$ such that $x^{T}Ax=0$, then $A\in PSD_{2,n}$. For higher order copositive tensors, we have the following property which is weaker than the one in the matrix case.\vskip 2mm

\begin{proposition}\label{copositive} If $\A \in COP_{m,n}$ and there exists some $x\in \R^n_{++}$ such that $\A x^m=0$, then $\A x^{m-1}=0$ and $\A x^{m-2}\in PSD_{2,n}$.
\end{proposition}
\noindent{\bf Proof.} For any given $u\in \R^n$, there exists some $\alpha_0>0$ such that for any $\alpha\in (0,\alpha_0)$, $x+\alpha u\in \R^n_+$. By the copositivity of $\A$, it follows readily that
$$ 0  \leq   \A (x+\alpha u)^m =\A x^m+ \sum\limits_{k=1}^m \alpha^k \A x^{m-k}u^k =\alpha \left(\A x^{m-1}u+\sum\limits_{k=2}^m \alpha^{k-1} \A x^{m-k}u^k \right),
$$
which further implies that
$$\A x^{m-1}u+\sum\limits_{k=2}^m \alpha^{k-1} \A x^{m-k}u^k \geq 0, \forall u\in \R^n. \eqno(5.3)$$
Let $\alpha\rightarrow 0$. (5.3) immediately yields that $\A x^{m-1}u\geq 0$ for all $u\in\R^n$. This indicates that $\A x^{m-1}=0$. Substituting this into (5.3), we get  $\A x^{m-2}u^2+\sum\limits_{k=3}^m \alpha^{k-2} \A x^{m-k}u^k \geq 0$, for all $u\in \R^n$. Thus $\A x^{m-2}$ is a positive semidefinite matrix. \qed

The above proposition also indicates that for a copositive tensor satisfying $\A x^m=0$ with some positive $x$, then $0$ is an $H^{++}$-eigenvalue of $\A$ and $x$ is the corresponding $H^{++}$-eigenvector.\vskip 2mm

The gap existing between doubly nonnegative matrices and completely positive matrices has been extensively studied \cite{BAD2009,BM2003,DA2013}. The remaining part of this section will be devoted to the equivalence and the gap between the tensor cones $DNN_{m,n}$ and $CP_{m,n}$. It is known from the literature of matrices that any rank-one matrix is completely positive if and only if it is nonnegative. This also holds for higher order tensors as the following proposition demonstrated.

\begin{proposition}\label{cp6} A rank-one symmetric tensor is completely positive if and only if it is nonnegative.
\end{proposition}
\noindent{\bf Proof.} The necessity is trivial by definition. To show the sufficiency, note that for any rank-one symmetric tensor $\A=\lambda x^m$ to be nonnegative, we have $x\neq 0$, $\lambda\neq 0$ and $\lambda x_{i_1}\cdots x_{i_m}\geq 0$, for all $i_1,\ldots, i_m\in[n]$. If $x$ has only one nonzero element, the desired statement holds immediately. If there exists at least two nonzero elements, we claim that all nonzero elements should be of the same sign. Otherwise, if $x_i>0$ and $x_j<0$, then $\lambda x_i^{m-1}x_j$ and $\lambda x_i^{m-2}x_j^2$ will not be nonnegative simultaneously. Thus all elements in $x$ are either nonnegative or nonpositive. When $m$ is even, we can easily get $\lambda>0$. Thus $\A$ is completely positive. If $m$ is odd, we can get that $\lambda^{1/m} x\geq 0$. Thus $\A$ is completely positive.\qed

The above proposition provides a special case that $CP_{m,n}$ coincides with $DNN_{m,n}$. Generally, there exists a gap between $DNN_{m,n}$ and $CP_{m,n}$. For example, if $\Q$ is the signless Laplacian tensor of a nonempty $m$-uniform hypergraph with $m\geq 3$, Theorem \ref{DNN-Tensors} and Proposition \ref{signless} lead to $\Q\in DNN_{m,n}\setminus CP_{m,n}$. Recall from \cite{BM2003} that for any matrix $A\in \S_{2,n}$, if $A$ is of rank $2$ or $n\leq 4$ , $A\in DNN_{2,n}$ if and only if $A\in CP_{2,m}$. In other words, $DNN_{2,n}=CP_{2,n}$ for these two cases. How about higher-order tensors? We answer this question in a negative way as follows.




\begin{proposition}\label{rank2} Let $m\geq 3$ and $n\geq 2$. Then
$$\left\{\alpha\left(e^{(i)}-e^{(j)}\right)^m+\alpha e^m: i,j\in [n], i\neq j, \alpha\in\R_{++}\right\}\subseteq DNN_{m,n}\setminus CP_{m,n}.$$
\end{proposition}
\noindent{\bf Proof.} For simplicity, denote $GAP:=\left\{\alpha\left(e^{(i)}-e^{(j)}\right)^m+\alpha e^m: i,j\in [n], i\neq j, \alpha\in\R_{++}\right\}$. Then for any $\A=\left(a_{i_1\ldots i_m}\right)\in GAP$, it follows from Proposition \ref{add} that $\A\in DNN_{m,n}$. Additionally, it is easy to verify that $\A$ is rank-two. However, as we can see, $a_{i\ldots ij}=0$ and $a_{i\ldots ijj}=2$. This indicates that $\A$ breaks the zero-entry dominance property. Thus $\A\in DNN_{m,n}\setminus CP_{m,n}$.\qed

\section{More Subclasses of DNN and CP Tensors}
More doubly nonnegative tensors and completely positive tensors are discussed in this section. We start with Cauchy tensors.

\begin{theorem}\label{Cauchy1}  Let $\C\in\S_{m,n}$ be a Cauchy tensor and $c=(c_1,\cdots, c_n)^{T}\in\R^n$ be its generating vector. The following statements are equivalent:

(i) $\C$ is completely positive;

(ii) $\C$ is strictly copositive;

(iii) $c>0$;

(iv) the function $f_{\C}(x):=\C x^m$ is strictly monotonically increasing in $\R^n_+$;

(v) $\C$ is doubly nonnegative.
\end{theorem}

\noindent{\bf Proof.} The implication ``$(ii)\Rightarrow (iii)$" follows readily from $0<\C e_i^m=\frac{1}{mc_i}$ for any $i\in [n]$. To get ``$(iii)\Rightarrow (i)$", we can employ the proof in Theorem 3.1 in \cite{CLQ2015} that for any $x\in \R^n$, it yields that

\begin{eqnarray}
   \C x^m &=& \sum\limits_{i_1,\cdots,i_m\in [n]} c_{i_1\cdots i_m}x_{i_1}\cdots x_{i_m} =\sum\limits_{i_1,\cdots,i_m\in [n]} \frac{x_{i_1}\cdots x_{i_m}}{c_{i_1}+\cdots+c_{i_m}} \nonumber \\
   &=& \sum\limits_{i_1,\cdots,i_m\in [n]}\int_{0}^{1} t^{c_{i_1}+\cdots+c_{i_m}-1} x_{i_1}\cdots x_{i_m} dt \nonumber  \\
  &=&\int_{0}^{1}\left( \sum\limits_{i_1,\cdots,i_m\in [n]} t^{c_{i_1}+\cdots+c_{i_m}-1} x_{i_1}\cdots x_{i_m}\right) dt \nonumber   \\
    &=& \int_{0}^{1}\left(\sum\limits_{i=1}^n t^{c_i-\frac{1}{m}}x_i\right)^m dt. \nonumber
\end{eqnarray}  Note that
\begin{eqnarray}
   \int_{0}^{1}\left(\sum\limits_{i=1}^n t^{c_i-\frac{1}{m}}x_i\right)^m dt &=& \lim\limits_{k\rightarrow \infty}\sum\limits_{j\in [k]} \left(\sum\limits_{i=1}^n\left(\frac{j}{k}\right)^{c_i-\frac{1}{m}}x_i\right)^m/k \nonumber \\
     &=&  \lim\limits_{k\rightarrow \infty}\sum\limits_{j\in [k]} \left(\sum\limits_{i=1}^n\left(\frac{j}{k}\right)^{c_i-\frac{1}{m}}x_i/k^{\frac{1}{m}}\right)^m  \nonumber  \\
    &=:& \lim\limits_{k\rightarrow \infty}\sum\limits_{j\in [k]}(\langle u^j,x\rangle)^m, \nonumber
\end{eqnarray}  with $$u^j:=\left(\frac{\left(j/k\right)^{c_1-\frac{1}{m}}}{k^{\frac{1}{m}}},\ldots, \frac{\left(j/k\right)^{c_n-\frac{1}{m}}}{k^{\frac{1}{m}}}\right)^T\in\R^n_+, \forall j\in [k].$$ By setting $\C_k:=\sum\limits_{j\in [k]}(u^j)^m$, it follows that $\C=\lim\limits_{k\rightarrow \infty} \C_k$ and $\C_k\in CP_{m,n}$. The closedness of $CP_{m,n}$ leads to $\C\in CP_{m,n}$. This implies (i). Conversely, if (i) holds, then $\C$ is certainly copositive, which deduces that $0\leq \C e_i^m=\frac{1}{mc_i}$, for all $i\in[n]$. Thus (iii) holds. Next we prove the equivalence between (iii) and (iv). Assume that (iii) holds, for any distinct $x$, $y\in\R^n_+$, satisfying $x\geq y$, i.e., there exists an index $i\in [n]$ such $x_i>y_i$, we have
$$  f_{\C}(x)-f_{\C}(y) ={\C}x^m-{\C}y^m=\sum\limits_{\begin{subarray}{c}
i_1,\cdots,i_m\in [n]\\
 (i_1,\cdots, i_m)\neq (i,\cdots,i)
 \end{subarray}
}\frac{x_{i_1}\cdots x_{i_m}-y_{i_1}\cdots y_{i_m}}{c_{i_1}+\cdots +c_{i_m}}+\frac{x_i^m-y_i^m}{mc_i}>0.$$
Thus (iv) is obtained. Conversely, if $f_{\C}(x)$ is strictly monotonically increasing in $\R^n_+$, then for any $i\in [n]$,
$0<f_{\C}(e_i)-f_{\C}(0)=\frac{1}{mc_i}$, which implies that $c>0$. Besides, by setting $x\in\R^n_+\backslash \{0\}$ and $y=0$, the strict monotonically increasing property of $f_{\C}$ also implies that $\C x^m>0$. Thus (iii) and (ii) hold. Trivially, we can get (v) from (i). If (v) holds, then $\C$ is copositive and hence for any $i\in [n]$, $0\geq \C e_i^m=\frac{1}{mc_i}$, which implies (iii). This completes the proof. \qed

\begin{proposition}\label{Cauchy2} For any given Cauchy tensor $\C\in\T_{m,n}$ with its generating vector $c=(c_1,\cdots, c_n)^{T}\in\R^n$. If $c>0$, then the following statements are equivalent:

(i) $c_1$, $\ldots$, $c_n$ are mutually distinct;

(ii) $\C$ is strongly doubly nonnegative;

(iii) $\C$ is strongly completely positive.
\end{proposition}
\noindent{\bf Proof.} It follows from Proposition \ref{Cauchy1} that $\C$ is completely positive and hence doubly nonnegative. When $m$ is even, the desired equivalence can be derived from Theorem 2.3 in \cite{CQ2015} and Proposition \ref{ppm}. Now we consider the case that $m$ is odd. To show the implication ``(i)$\Rightarrow$ (ii)", we assume on the contrary that $0$ is an $H$-eigenvalue of $\C$ with its associated $H$-eigenvector $x$. Then for any $i\in [n]$, we have
 \begin{eqnarray}
   0=\left(\C x^{m-1}\right)_i &=& \sum\limits_{i_2,\ldots,i_m\in[n]} \frac{x_{i_2}\cdots x_{i_m}}{c_i+c_{i_2}+\cdots+c_{i_m}} \nonumber \\
    &=& \sum\limits_{i_2,\ldots,i_m\in[n]}  \int_0^1 t^{c_i+c_{i_2}+\cdots+c_{i_m}-1} x_{i_2}\cdots x_{i_m}  dt \nonumber\\
    &=& \int_0^1 t^{c_i}\left(\sum\limits_{j\in [n]}t^{c_j-\frac{1}{m-1}}x_j\right)^{m-1} dt, \nonumber
 \end{eqnarray}
which implies that $\sum\limits_{j\in [n]}t^{c_j-\frac{1}{m-1}}x_j\equiv 0$, for all $t\in [0,1]$. Thus,
$$x_1+t^{c_2-c_1}x_2+\cdots +t^{c_n-c_1}x_n=0,~\forall t\in [0,1].$$
By the continuity and the condition that all components of $c$ are mutually distinct, it follows readily that $x_1=0$. Then we have $x_2+t^{c_3-c_2}x_2+\cdots +t^{c_n-c_2}x_n=0$, $\forall t\in [0,1],$ which implies $x_2=0$. By repeating this process, we can gradually get $x=0$, which contradicts to the assumption that $x$ is an $H$-eigenvalue. Thus (ii) is obtained. Conversely, to show ``(ii)$\Rightarrow$ (i)", we still assume by contrary that $c_1$, $\ldots$, $c_n$ are not mutually distinct. Without loss of generality, we assume that $c_1=c_2$. By setting $x\in \R^n$ with $x_1=-x_2=1$ and other components $0$, we find that for any $i\in [n]$,
$$
  \left(\C x^{m-1}\right)_i= \int_0^1 t^{c_i}\left(\sum\limits_{j\in [2]}t^{c_j-\frac{1}{m-1}}x_j\right)^{m-1} dt=
\int_0^1 t^{c_i}\left(t^{c_1-\frac{1}{m-1}}-t^{c_2-\frac{1}{m-1}}\right)^{m-1} dt= 0,$$
which indicates that $0$ is an $H$-eigenvalue of $\C$. This is a contradiction to the condition that $\C\in SDNN_{m,n}$. Thus the desired implication holds. Note that $\C$ is completely positive since $c>0$, the equivalence of (ii) and (iii) can be achieved by applying Theorem \ref{SDNN-SCP}. \qed

As said before, Fan and Zhou \cite{FZ2014} proposed an optimization algorithm to decompose a completely positive tensor. This is an important contribution to the study of completely positive tensor decomposition. Now, our work shows that positive Cauchy tensors can serve as testing examples for completely positive tensor decomposition.   This is confirmed by computation \cite{F2015}.

Recall that a tensor $\A=\left(a_{i_1\ldots i_m}\right)\in\S_{m,n}$ is called a Hilbert tensor if $a_{i_1\ldots i_m}=\frac{1}{ i_1+\cdots +i_m-m+1 }$ (\cite{SongQi2014}). Obviously, a Hilbert tensor is both a Cauchy tensor and a Hankel tensor \cite{CQ2015}.

\begin{corollary}\label{Hilbert} A Hilbert tensor is strongly completely positive.
\end{corollary}

Nonnegative strong Hankel tensors are doubly nonnegative but not always completely positive. Several sufficient conditions to ensure the complete positivity of a nonnegative strong Hankel tensor were proposed in \cite{DQW2015,LQX2014} as recalled below.\vskip 2mm

\begin{theorem}\label{Hankel1} Let $\A=(a_{i_1\ldots i_m})\in \S_{m,n}$ be a nonnegative Hankel tensor with its generating vector $h=(h_0,\ldots, h_{m(n-1)})^T$ and its generating function $f$. Then $\A\in CP_{m,n}$ if one of the following holds

(i) $f$ is nonnegative, and $f(t)=0$ for any $t<0$;

(ii) if the $m$th order $(2n-1)$-dimensional Hankel tensor generated by $\widetilde{\bf h} = (h_0,0,h_1,0,h_2,\dots,0,h_{m(n-1)})^T$ is a strong Hankel tensor.

\end{theorem}

Some necessary conditions for completely positive Hankel tensors are presented as below.

\begin{proposition}\label{Hankel2} Let $\A\in \S_{m,n}$ be a Hankel tensor. If $\A$ is completely positive, then all its induced tensors are completely positive.
\end{proposition}
\noindent{\bf Proof.} It is known from Theorem 4.1 in \cite{Q2015} that every Hankel tensor has a Vandermonde decomposition. Thus we write $\A$ as
$\A=\sum\limits_{k=1}^r\alpha_k \left(u^{(k)}\right)^m$, where $\alpha_k\in \R\setminus\{0\}$, $u^{(k)}=\left(1,\mu_k,\ldots, \mu_k^{n-1}\right)^{T}\in\R^n$, $k\in [r]$. For any $s\in [m]$, the corresponding $s$-induced tensor $\B_s=\sum\limits_{k=1}^r\alpha_k \left(u^{(k)}\right)^s$ can be rewritten as
$\B_s=\sum\limits_{k=1}^r\alpha_k \left((e^{(1)})^{T}u^{(k)}\right)^{m-s}\left(u^{(k)}\right)^s$. Thus, the desired result follows from (ii) of Proposition \ref{cp3} by setting $x=e^{(1)}$. \qed

\section{Applications}
\subsection{Application 1: Preprocessing for CP Tensors}
The completely positive tensor verification and decomposition are very important as discussed in \cite{FZ2014,K2015}. In this section, by employing the zero-entry dominance property and a simplified strong dominance property called the one-duplicated dominance property, a preprocessing scheme is proposed to accelerate the verification procedure for completely positive tensors based on the Fan-Zhou algorithm. The one-duplicated property is defined as follows.

\begin{definition}\label{One-duplicated} Let $i_j \in [n]$ for $j \in [m]$.  We say that $(l_1, \cdots, l_m)$ is one-duplicated from $(i_1, \cdots, i_m)$, if $l_j = i_j$ for $j \in [m], j \not = p$ and $l_p = i_q$ for $q \in [m], q \not = p$.  Denote $S_1(i_1, \cdots i_m)$ be the set of $(l_1, \cdots, l_m)$, where $(l_1, \cdots, l_m)$ is one-duplicated from $(i_1, \cdots, i_m)$. A nonnegative tensor $\A\in \S_{m,n}$ is said to have the \emph{one-duplicated dominance property} if for any $(i_1, \ldots, i_m)$, where $i_j \in [n]$ for $j \in [m]$,
 $$a_{i_1\ldots i_m} \le {1 \over m(m-1)} \sum \left\{ a_{j_1\ldots, j_m} : (j_1, \cdots, j_m) \in S_1(i_1, \cdots, i_m) \right\}.\eqno(7.1)$$
\end{definition}
Easily, one-duplicated dominance property can be derived from the strong dominance property with $s=m(m-1)$ in Proposition \ref{strong dominance property}, and hence it is a necessary condition for tensors to be completely positive.

The preprocessed Fan-Zhou algorithm to verify a given nonnegative symmetric tensor $\A=\left(a_{i_1\ldots i_m}\right)\in \S_{m,n}$ is presented as follows.

\begin{algorithm}[htbp]
\caption{The Preprocessed Fan-Zhou Algorithm.}
\label{alg_vander}
\begin{algorithmic}[1]
\Require
A nonnegative symmetric tensor $\A$;
\Ensure
Certificate for the Non-Complete-Positivity or a CP-tensor decomposition of $\A$;

\begin{itemize}
\item[Step~0] Set $\epsilon >0$. Denote $O(\A):=\{(i_1,\ldots,i_m): a_{i_1\ldots i_m}=0\}$. If $O(\A)=\emptyset$, go to Step 1. Otherwise, check the zero-entry dominance property: Let $\B=\left(b_{i_1,\ldots,i_m}\right)=\A$. Then, pick any $(i_1,\ldots,i_m)\in O(\A)$, and set $b_{j_1,\ldots,j_m}=0$ for all $\{j_1,\ldots,j_m\}\supseteq\{i_1,\ldots,i_m\}$. If $\|\A-\B\|_F>\epsilon$, then $\A\notin CP_{m,n}$ and stop; otherwise go to Step 1.

\item[Step~1] Check the one-duplicated dominance (7.1). If not satisfied, then $\A\notin CP_{m,n}$ and stop; otherwise go to Step 2.

\item[Step~2] Use the Fan-Zhou algorithm in \cite{FZ2014}.

\end{itemize}

\end{algorithmic}
\end{algorithm}

%
%
%

%
%
%

Some numerical tests are reported as follows.

\begin{example}\label{Strong Hankel}{\bf[Hankel Tensors]} Randomly generate a vector $\xi=(\xi_1,\ldots, \xi_r)^T\in\R^r$ and use it to get a complete Hankel tensor $\B(\xi)=\left(b_{i_1,\ldots,i_m}(\xi)\right)\sum\limits_{i\in [r]}\left(u^{(k)}\right)^m$ with $u^{(k)}=(1,\xi,\ldots, \xi^{n-1})^T$. Let $\alpha(\xi)=\min\{t(\xi) ,0\}$ with $t(\xi)$ the minimal entry of $\B(\xi)$,  and set $\A(\xi)=\B(\xi)-\alpha(\xi) e^m$. Such a tensor $\A(\xi)$ is always a nonnegative complete Hankel tensor and hence doubly nonnegative by Theorem \ref{DNN-Tensors}. By applying Algorithm 1 with $\epsilon=1{\rm e}\hspace{-2pt}-\hspace{-2pt}12$, we find that the majority of such tensors $\A(\xi)$'s can be excluded from the class of completely positive tensors by Step $0$ and Step $1$ as the following table illustrated.

\begin{table}[h]
\begin{center}\caption{ Preprocessing for Hankel tensors}
\begin{tabular}{c|c|c|c|c}
   \hline
   \hline
    No. &  m &  n & r  &  Perc. \\ \hline
   1,000,000  &  3 & 3  &  3 &  79.0\%  \\  
   100,000  &  5 & 3  &  6 &  90.2\% \\  
    100,000 &  11 & 3  & 6  &  91.0\% \\   
    10,000 &  12 & 4 & 7 &  92.4\% \\   
    10,000 &  3 & 11 &  16 &  58.2\% \\   
    10,000 &  4 &  12 &  18 &  56.1\% \\   
   \hline
   \hline
 \end{tabular}

\end{center}
\end{table}
\noindent Here ``No." stands for the number of randomly generated tensors and ``Perc." presents the percentage that can be excluded from Step 0 and Step 1 when applying Algorithm 1. This indicates that our preprocessing scheme is very efficient and will greatly accelerate Fan-Zhou algorithm for handling with complete Hankel tensors. Besides, this table also tells us that most of the nonnegative complete Hankel tensors lie in the gap between $DNN_{m,n}$ and $CP_{m,n}$, with some special cases which are completely positive. This phenomenon is partially illustrated with the following selected cases with $m=3$, $n=11$ and $r=16$.

\begin{table}[h]
\begin{center}\caption{ Selected results for Hankel tensors }
\begin{tabular}{c|c|c|c|c}
    \hline
   \hline
     $\xi$ & $\alpha(\xi)$ \ &CP & Not CP & Excluded by\\ \hline
         $(1.3769,   -1.4082,    0.1412,    1.2897,   -0.4949,   -0.6248,$ &   & &   & \\
    $   -0.9336,   -0.2787, -0.2005,    0.1367,   -0.8833,$   &  $-8.2003e+03$ & & $\surd$  & Step $0$ \\
    $    0.0825,   -0.6039,   -0.3687,   -0.8382,   -0.2825)^T$ &   &  & &\\
   \hline
         $(-0.7841,   -1.8054,    1.8586,   -0.6045,    0.1034,    0.5632,$ &   & &   &  \\
    $ 0.1136,   -0.9047, -0.4677,   -0.1249,    1.4790,$   &  $-1.6321$ & & $\surd$  & Step $1$ \\
    $ -0.8608,    0.7847,    0.3086,   -0.2339,   -1.0570)^T$ &   &  & &\\
   \hline
         $(2.5610,   0.1966,    0.7577,    2.0048,   0.9201,    1.6254,$ &   & &   &  \\
    $ 1.7530,    1.2135,  0.2298,    0.9929,    1.0932,$   &  $0$ & $\surd$ &  &   \\
    $ 1.9353,    1.6635,    0.6498,    2.6199,    0.9492)^T$ &   &  & &\\
   \hline
    \hline
 \end{tabular}

\end{center}
\end{table}
\end{example}

\begin{example}\label{SNT}{\bf[Symmetric Nonnegative Tensors]} Note that $$BASIS:=\left\{\left(\sum\limits_{j\in [m]}e^{(i_j)}\right)^m: i_1,\ldots,i_m\in [n], i_1\leq \cdots \leq i_m\right\}$$ is a basis of $\S_{m,n}$. That means, any tensor $\B\in\S_{m,n}$ can be written as a linear combination of elements in $BASIS$. Henceforth, we randomly generate a symmetric nonnegative tensor as follows:
$$\B(\alpha)=\sum\limits_{i_1,\ldots,i_m\in [n], i_1\leq \cdots \leq i_m} \alpha_{i_1\ldots i_m}\left(\sum\limits_{j\in [m]}e^{(i_j)}\right)^m,$$
where $\alpha_{i_1\ldots i_m}=randn+t$ with a random scalar $randn$ obeying the standard normal distribution for all $i_1,\ldots,i_m\in [n]$, and some positive scalar $t$. Similar as in Example \ref{Strong Hankel}, let $\rho(\alpha)=\min\{\gamma(\alpha) ,0\}$ with $\gamma(\alpha)$ be minimal entry of $\B(\alpha)$,  and set $\A(\alpha)=\B(\alpha)-\rho(\alpha) e^m$. Such a tensor $\A(\alpha)$ is symmetric and nonnegative. Set different values to $t$ and randomly generate {\bf 10,000} $\A$'s with each $t$. By applying Algorithm 1, we find that a considerable percentage of tensors can be efficiently excluded by Step 0 and Step 1, as the following table shows.
\begin{table}[h]
\begin{center}\caption{ Preprocessing for symmetric nonnegative tensors}
\begin{tabular}{c|r|r|r|r|r|r|r|r}
   \hline
   \hline
    m  & 5       & 5        & 5         &  11      & 11         & 12           & 3      & 4     \\
    n  & 3       & 3        & 3         &  3       & 3          & 4            & 11      & 10    \\
    t  & 0.1     & 0.4      & 1         &  0.4     & 1          & 0.4          & 0.4      &  0.4     \\     \hline
Perc.  &99.5\%   & 93.8\%   & 40.7\%    & 96.7\%   & 24.5\%     & 85.1\%             & 22.4\%   & 37.7\%    \\
       \hline
       \hline
 \end{tabular}
\end{center}
\end{table}
\end{example}

These examples show that our preprocessing scheme may accelerate the verification for completely positive tensors based on the Fan-Zhou algorithm efficiently.

\subsection{Application 2: Tensor Complementarity Problems}

Recently, the tensor complementarity problem (TCP) is studied \cite{CQW, LQX, SQ2014, SQ2015}. Let $\A \in \T_{m, n}$ and $q \in \R^n$. The tensor complementarity problem TCP($q, \A$) is to find $x \in \R^n$ such that
$$x \ge 0, \ q + \A x^{m-1} \ge 0, \ {\rm and} \ x^\top(q + \A x^{m-1}) = 0.$$

For strongly doubly nonnegative tensors, we have the following property on the corresponding tensor complementarity problem.


\begin{proposition} Let $\A\in \S_{m,n}$ and $q\in\R^n$. If $\A\in SDNN_{m,n}$, then TCP($q, \A$) has a nonempty, compact solution set.
\end{proposition}

\noindent{\bf Proof.} It follows readily from Theorem 5 in \cite{CQW} and (ii) of Proposition \ref{cone1}. \qed

\section{Conclusions}
In this paper, the double nonnegativity is extended from matrices to tensors of any order in terms of the nonnegativity of entries and $H$-eigenvalues, and the completely positive tensors,
as a very important subclass of doubly nonnegative tensors, are further studied. Our contributions are three-fold. Firstly, several structured tensors, which have wide applications in many
real-life problems, are shown to be doubly nonnegative in both even and odd order cases. It is worth pointing out that for odd order tensors, the proposed double nonnegativity, to some extent,
can be served as a counterpart of the positive semidefiniteness property for even order tensors. This makes up the deficiency that the positive semidefiniteness property vanishes in the odd
order case. Secondly, for completely positive tensors, the dominance properties are exploited to exclude a number of symmetric nonnegative tensors, such as the signless Laplacian tensors
of nonempty $m$-uniform hypergraphs with $m\geq 3$, from the class of completely positive tensors. Moreover, these properties are also used in our preprocessing scheme to accelerate the
verification for completely positive tensors based on the Fan-Zhou algorithm. Thirdly, all positive Cauchy tensors of any order (even or odd) are shown to be completely positive,
which serves as a new sufficient condition and provides an easily verifiable structure in the study of completely positive tensors and decomposition. In addition, the solution analysis for tensor complementarity problems with the strongly doubly nonnegative tensor structure is discussed. All these results can be served as a supplement to enrich tensor analysis, computation and applications.

\vskip 4mm

\section*{Acknowledgements}
The authors would like to thank Dr. Yannan Chen for his support on the numerical tests, Prof. Yiju Wang and Mr. Weiyang Ding for their valuable comments. We appreciate Prof. Masakazu Kojima and Prof. Kim-Chuan Toh for sharing their current work on polynomial optimization problems. We thank Prof. Jinyan Fan for the discussion and her numerical experiments on positive Cauchy tensors.


\end{document}